\documentclass[12pt]{article}
\usepackage{amssymb,graphicx}
\def\no{\noindent}

\def\f#1#2{\frac{#1}{#2}}
\def\f{\frac}

\def\ov{\overline}

\def\dint{\int\!\!\int}

\def\R{{\mathbb R}}

\def\C{{\cal C}}
\def\E{{\cal E}}

\def\p{\partial}

\def\vp{\varphi}
\def\L{{\bf L}}

\def\O{{\cal O}}

\def\v{\vskip 1em}
\def\vs{\vskip 2em}
\def\vsk{\vskip 4em}

\def\ve{\varepsilon}

\def\meas{\hbox{meas}\,}

\def\T{{\cal T}}

\def\sqr#1#2{\vbox{\hrule height .#2pt
\hbox{\vrule width .#2pt height #1pt \kern #1pt
\vrule width .#2pt}\hrule height .#2pt }}
\def\square{\sqr74}
\def\endproof{\hphantom{MM}\hfill\llap{$\square$}\goodbreak}

\newcommand{\beq}{\begin{equation}}
\newcommand{\eeq}{\end{equation}}
\newcommand{\ben}{\begin{eqnarray}}
\newcommand{\een}{\end{eqnarray}}
\newcommand{\beno}{\begin{eqnarray*}}
\newcommand{\eeno}{\end{eqnarray*}}

\begin{document}
\title{\bf  Global Conservative Solutions to a Nonlinear 
Variational Wave Equation}
\author{Alberto Bressan and Yuxi Zheng \\[2mm]
{\small Department of Mathematics, The Pennsylvania State University }\\[-2mm]
{\small E-mail: bressan@math.psu.edu; yzheng@math.psu.edu} }
\maketitle

\begin{abstract}
We establish the existence of a conservative weak solution to the
Cauchy problem for the nonlinear variational wave equation 
$u_{tt} - c(u)(c(u)u_x)_x=0$, for initial data of finite energy. 
Here  $c(\cdot)$ is any smooth function with uniformly positive
bounded values. 
\end{abstract}

\vspace*{3mm}


\renewcommand{\theequation}{\thesection.\arabic{equation}}

\section{Introduction}

\setcounter{equation}{0}

We are interested in the Cauchy problem
\beq
u_{tt} - c(u)\big(c(u) u_x\big)_x=0\,,\label{1.1}
\eeq
with initial data
\beq
u(0,x)=u_0(x)\,,\qquad
u_t(0,x)=u_1(x)\,.\label{1.2}
\eeq
Throughout the following, we assume that $c:\R\mapsto \R_+$ is
a smooth, bounded, uniformly positive function.
Even for smooth initial data, it is well known that the solution
can lose regularity in finite time (\cite{[GHZ]}).  It is thus of interest
to study whether the solution can be extended beyond the
time when a singularity appears.  This is indeed the main
concern of the present paper.  In (\cite{[BZZ]}) we considered the
related equation
\beq
u_t+f(u)_x={1\over 2} \int_0^x f''(u)u_x^2\,dx\label{1.3'}
\eeq
and constructed a semigroup of solutions, depending continuously
on the initial data.
Here we establish similar results for the nonlinear wave equation
(\ref{1.1}).  By introducing new sets of dependent and independent variables,
we show that the solution to the Cauchy problem can be
obtained as the fixed point of a contractive transformation.
Our main result can be stated as follows.
\v
\no{\bf Theorem 1.} {\it Let $c:\R\mapsto [\kappa^{-1},\,\kappa]$ be
a smooth function, for some $\kappa>1$.  Assume that the initial data
$u_0$ in (\ref{1.2}) is absolutely continuous, and that
$(u_0)_x\in\L^2\,$,   ~$u_1\in\L^2$.
Then the Cauchy problem (\ref{1.1})-(\ref{1.2}) admits a weak solution
$u=u(t,x)$, defined for all $(t,x)\in\R\times\R$.
In the $t$-$x$ plane, the function $u$ is locally H\"older continuous
with exponent $1/2$.   This solution $t\mapsto u(t,\cdot)$ is
continuously differentiable as a map with values in $\L^p_{\rm loc}$, for all
$1\leq p<2$. Moreover, it is Lipschitz continuous w.r.t.~the $\L^2$
distance, i.e.
\beq
\big\|u(t,\cdot)-u(s,\cdot)\big\|_{\L^2}\leq L\,|t-s|\label{1.lip}
\eeq
for all $t,s\in\R$.
The equation (\ref{1.1}) is satisfied in integral sense, i.e.
\beq
\dint \Big[\phi_t\, u_t - \big(c(u) \phi\big)_x c(u)\,u_x\Big]\,dxdt =0
\label{1.4'}
\eeq
for all test functions $\phi\in\C^1_c$.  Concerning the initial conditions, 
the first equality in (\ref{1.2})
is satisfied pointwise, while the second holds in  $\L^p_{\rm loc}\,$
for $p\in [1,2[\,$.
}
\v
Our constructive procedure yields solutions which depend continuously
on the initial data.  Moreover, the ``energy"
\beq
\E(t)\doteq \f12\int \Big[u_t^2(t,x) +c^2\big(u(t,x)\big) \,u_x^2(t,x)\Big]
\,dx\label{1.5'}
\eeq
remains uniformly bounded. More precisely, one has
\v
\no{\bf Theorem 2.} {\it A family of weak solutions to the
Cauchy problem (\ref{1.1})-(\ref{1.2}) can be constructed with the
following additional properties.
For every $t\in\R$ one has
\beq
\E(t)~\leq ~\E_0~\doteq~ \f12\int \Big[ u_1^2(x)+c^2\big(u_0(x)\big) 
(u_0)_x^2(x) \Big]\,dx\,.\label{1.7}
\eeq
Moreover, let a sequence of initial conditions satisfy
$$
\big\|(u_0^n)_x- (u_0)_x\big\|_{\L^2}\to 0\,,\qquad\qquad
\big\|u_1^n- u_1\big\|_{\L^2}\to 0\,,
$$
and $u^n_0\to u_0$ uniformly on compact sets, as $n\to \infty$.
Then one has the convergence of the corresponding solutions $u^n\to u$,
uniformly on bounded subsets of the $t$-$x$ plane.}
\v
It appears in (\ref{1.7}) that the total energy of our solutions
may decrease in time.  Yet, we emphasize that
our solutions are {\it conservative}, in the following sense.
\v
\no{\bf Theorem 3.} {\it There exists a continuous family
$\{\mu_t\,;~~t\in\R\}$ of
positive Radon measures on the real line with the following properties.
\begin{description}
\item[](i) At every time $t$, one has $\mu_t(\R)=\E_0$.

\item[](ii) For each $t$, the absolutely continuous part of $\mu_t$
has density $\f12(u_t^2+c^2\,u_x^2)$ w.r.t.~the Lebesgue measure.

\item[](iii)  For almost every $t\in\R$, the singular part of $\mu_t$
is concentrated on the set where $c'(u)=0$.
\end{description}
}
In other words, the total energy represented by the measure $\mu$ is
conserved in time.  Occasionally, some of this energy is concentrated
on a set of measure zero. At the times $\tau$ when this happens,
$\mu_\tau$ has a non-trivial singular part and
$\E(\tau)<\E_0$.  The condition (iii) puts some restrictions
on the set of such times $\tau$. In particular, if $c'(u)\not= 0$
for all $u$, then this set has measure zero.
\v
The paper is organized as follows.  In the next two subsections
we briefly
discuss the physical motivations for
the equation and recall some known results on its solutions.
In Section 2 we introduce a new set of independent and dependent
variables, and derive some identities
valid for smooth solutions.  We formulate a set of equations in the 
new variables which is equivalent to (\ref{1.1}).  Remarkably,
in the new variables all singularities disappear: Smooth 
initial data lead to
globally smooth solutions.
In Section 3 we use a contractive transformation in a Banach space 
with a suitable weighted norm to show that there is a unique
solution to the set of equations in the new variables,
depending continuously on the data $u_0, u_1$.
Going back to the original variables $u,t,x$,
in Section 4 we establish the H\"older continuity
of these solutions $u=u(t,x)$, and show
that the integral equation (\ref{1.4'}) is satisfied.
Moreover, in Section 5, we study the conservativeness of the solutions,
establish the energy inequality and the Lipschitz continuity of the map
$t\mapsto u(t,\cdot)$.  This already yields a proof
of Theorem 2. In Section 6 we study
the continuity of the maps
$t\mapsto u_x(t,\cdot)$, $t\mapsto u_t(t,\cdot)$, completing the proof
of Theorem 1.  The proof of Theorem 3 is given in  Section 7.

\subsection{Physical background of the equation}

Equation (\ref{1.1}) has several physical origins. 
In the context of nematic liquid crystals, 
it comes as follows. 
The mean orientation of the long molecules in a nematic liquid crystal
is described by a director field of unit vectors,  
${\mathbf n}\in{\mathbb S}^2$,
the unit sphere. Associated with the director field ${\mathbf n}$, there is
the well-known Oseen-Franck potential energy density $W$ given by
\beq
W\left({\mathbf n},\nabla{\mathbf n}\right) =
\alpha\left|{\mathbf n}\times(\nabla\times{\mathbf n})\right|^2
+ \beta(\nabla\cdot{\mathbf n})^2 +
\gamma\left({\mathbf n}\cdot\nabla\times{\mathbf n}\right)^2.
\label{2b}
\eeq
The positive constants $\alpha$, $\beta$, and $\gamma$ are elastic
constants of the liquid crystal. For the special case 
$\alpha = \beta = \gamma$, the potential energy density reduces to
$$
W\left({\mathbf n},\nabla{\mathbf n}\right) = \alpha\left|
\nabla{\mathbf n}\right|^2,
$$
which is the potential energy density used in harmonic maps into the sphere 
${\mathbb S}^2$.
There are many studies on the constrained elliptic system of equations for
${\mathbf n}$ derived through variational principles from the potential 
(\ref{2b}), and on the parabolic flow associated
with it, see \cite{[3],[7],[10],[19],[26],[49]} and references therein.
In the regime in which inertia effects dominate
viscosity, however, the propagation of the orientation waves in the
director field may then be modeled by the least action principle (Saxton 
\cite{[37]})
\beq
\frac{\delta}{\delta u} \int \left\{ \p_t{\mathbf n}\cdot\p_t{\mathbf n} 
- W\left({\mathbf n},
\nabla{\mathbf n}\right)\right\}\,d{\mathbf x}dt  = 0,
\qquad {\mathbf n}\cdot {\mathbf n} = 1.
\label{2a}
\eeq
In the special case $\alpha = \beta = \gamma$, this variational principle 
(\ref{2a}) yields the equation for harmonic wave maps from $(1+3)$-dimensional
 Minkowski space into the two sphere, see \cite{[6],[40],[41]} for example.
For planar deformations depending on a single space variable $x$,
the director field has the special form
$$
{\mathbf n} =  \cos u(x,t){\mathbf e}_x + \sin u(x,t){\mathbf e}_y,
$$
where the dependent variable $u\in\R^1$ measures the angle
of the director field to the $x$-direction, and ${\mathbf e}_x$ and
${\mathbf e}_y$ are the coordinate vectors in the $x$ and $y$ directions,
respectively.
In this case, the variational principle (\ref{2a}) reduces to (\ref{1.1})
with the wave speed $c$ given specifically by
\beq
c^2(u) = \alpha\cos^2u + \beta\sin^2u.
\label{2c}
\eeq

The equation (\ref{1.1}) has interesting connections with 
long waves on a dipole chain in the continuum limit (\cite{[16]},  
Zorski and Infeld \cite{[Zo]}, and Grundland and Infeld \cite{[17]}),
and in classical field theories and general relativity (\cite{[16]}).
We refer the interested reader to the article \cite{[16]} for these
connections.

This equation (\ref{1.1}) compares interestingly with other well-known
equations, e. g.
\beq
\p^2_tu -\p_x[p(\p_xu)] =0, \label{conlaw}
\eeq
where $p(\cdot)$ is a given function, considered by Lax \cite{[29]},
Klainerman and Majda \cite{[28]}, and Liu \cite{[34]}. Second related equation
is 
\beq
\p_t^2u- c^2(u) \Delta u = 0\label{lin}
\eeq
considered by  Lindblad \cite{[31]}, who established the global existence of 
smooth solutions of (\ref{lin}) with smooth, small, and spherically 
symmetric initial data in $\R^3$,
where the large-time decay of solutions in high space dimensions is crucial.
The multi-dimensional generalization of equation (\ref{1.1}),
\beq
\p_t^2u- c(u) \nabla\cdot \left(c(u) \nabla u\right) = 0,\label{multid}
\eeq
contains a lower order term proportional to $cc'|\nabla u|^2$,
which (\ref{lin}) lacks. This lower order term is responsible for
the blow-up in the derivatives of $u$. Finally, we note that equation 
(\ref{1.1}) also looks related to the perturbed wave equation
\beq
\p_t^2u - \Delta u+ f(u,\nabla u, \nabla \nabla u) = 0, \label{1.7cpde}
\eeq
where $f(u, \nabla u,\nabla\nabla u)$ satisfies an appropriate convexity
condition (for example, $f =u^p$ or $f=a (\p_tu)^2 + b|\nabla u|^2$) or
some nullity condition. Blow-up for (\ref{1.7cpde}) with a convexity condition
has been studied extensively, see 
\cite{[2],[14],[18],[23],[25],[30],[38],[43],[44]} and Strauss \cite{[45]}
for more reference. Global existence and uniqueness of solutions to 
(\ref{1.7cpde})
with a nullity condition depend on the nullity structure and large time decay
of solutions of the linear wave equation in higher dimensions
(see Klainerman and Machedon \cite{[27]} and references therein).
Therefore (\ref{1.1}) with the dependence of $c(u)$ on $u$ and the possibility
of sign changes in $c'(u)$ is familiar yet truly different.

Equation (\ref{1.1}) has interesting asymptotic uni-directional wave equations.
Hunter and Saxton (\cite{[21]}) derived the asymptotic equations
\beq
(u_t +u^nu_x)_x = \f12 nu^{n-1}(u_x)^2
\label{3cpde}
\eeq
for (\ref{1.1}) via weakly nonlinear geometric optics.
We mention that the $x$-derivative of equation
(\ref{3cpde}) appears in the high-frequency limit of the variational principle
for the Camassa-Holm equation (\cite{[1],[4],[5]}), which arises in the theory
of shallow water waves. A construction of global solutions to the
Camassa-Holm equations, based on a similar variable transformation 
as in the present paper, will appear in \cite{BrCo}

\subsection{Known results}
In \cite{[22]}, Hunter and Zheng established the global existence of weak
solutions to (\ref{3cpde}) ($n=1$) with initial data of bounded variations.
It has also been shown that the dissipative solutions are limits of
vanishing viscosity. Equation  (\ref{3cpde}) ($n=1$) is also shown
to be completely integrable (\cite{[HZ]}). 
In \cite{[52]}--\cite{[59]}, Ping Zhang and Zheng
study the global existence, uniqueness, and regularity of the weak
solutions to (\ref{3cpde}) ($n=1, 2$) with $L^2$ initial data, and
special cases of (\ref{1.1}). The study of
the asymptotic equation has been very beneficial for both the blow-up
result \cite{[GHZ]} and the current global existence result for the wave 
equation (\ref{1.1}).

\section{Variable Transformations} 

\setcounter{equation}{0}

We start by deriving some identities valid for smooth solutions.
Consider the variables
\beq
\left\{ 
\begin{array}{rcl}
R & \doteq  &u_t+c(u)u_x\,, \\
S & \doteq  &u_t-c(u)u_x\,,
\end{array} \right.\label{2.1}
\eeq
so that
\beq
u_t={R+S\over 2}\,,\qquad\qquad u_x={R-S\over 2c}\,.\label{2.2}
\eeq
By (\ref{1.1}), the variables $R,S$ satisfy
\beq
\left\{
\begin{array}{rcl}
R_t-cR_x &= & {c'\over 4c}(R^2-S^2), \\ [3mm]
S_t+cS_x &= & {c'\over 4c}(S^2-R^2).
\end{array} \right.
\label{2.3}
\eeq
Multiplying the first equation in (\ref{2.3}) by $R$ and the
second one by $S$, we obtain balance laws for $R^2$ and $S^2$, namely
\beq\left\{
\begin{array}{rcl}
(R^2)_t - (cR^2)_x & = & {c'\over 2c}(R^2S - RS^2)\, , \\ [3mm]
(S^2)_t + (cS^2)_x & = & - {c'\over 2c}(R^2S -RS^2)\,.
\end{array} 
\right.\label{2.4}
\eeq
As a consequence, the following quantities are conserved:
\beq
E\doteq {1\over 2}\big(u_t^2+c^2u_x^2\big)={R^2+S^2\over 4}\,,\qquad\qquad
M\doteq -u_tu_x={S^2-R^2\over 4c}\,.\label{2.5} 
\eeq
Indeed we have
\beq
\left\{ 
\begin{array}{rcl}
E_t+(c^2M)_x & = & 0\,, \\ 
M_t+E_x      & = & 0\,.
\end{array}\right.\label{2.6} 
\eeq
One can think of $R^2/4$ as the energy density of backward
moving waves, and $S^2/4$ as the energy density of forward moving waves.

We observe that, if $R,S$ satisfy (\ref{2.3}) and $u$ satisfies
(\ref{2.1}b), then the quantity
\beq
F \doteq  R - S -2cu_x \label{F}
\eeq
provides solutions to the linear homogeneous equation
\beq
F_t - c\,F_x = {c'\over 2c}(R+S+2cu_x)F. \label{FE}
\eeq
In particular, if $F\equiv 0$ at time $t=0$, the same
holds for all $t>0$.  Similarly, if $R,S$ satisfy (\ref{2.3}) and
$u$ satisfies (\ref{2.1}a), then the quantity
$$
G\doteq R+S-2u_t  
$$
provides solutions to the linear homogeneous equation
$$
G_t+c\,G_x={c'\over 2c}(R+S-2cu_x)G\,. 
$$
In particular, if $G\equiv 0$ at time $t=0$, the same
holds for all $t>0$.  We thus have
\v
\no{\bf Proposition 1.}
{\it Any smooth solution of (\ref{1.1}) provides a solution to
(\ref{2.1})--(\ref{2.3}). Conversely, any smooth solution of (\ref{2.1}b)
 and (\ref{2.3}) (or (\ref{2.1}a) and (\ref{2.3})) which satisfies
(\ref{2.2}b)  (or (\ref{2.2}a) ) at time $t=0$ provides
a solution to (\ref{1.1}). }
\vs

The main difficulty in the analysis of (\ref{1.1}) is the possible
breakdown of regularity of solutions.   Indeed, even for smooth
initial data, the quantities $u_x, u_t$ can blow up in finite time.
This is clear from the equations (\ref{2.3}), where the right
hand side grows quadratically, see (\cite{[GHZ]}) for handling change of 
signs of $c'$ and interaction between $R$ and $S$.
To deal with possibly unbounded values of $R,S$,
it is convenient to introduce a new set of dependent variables:
$$
w\doteq 2\arctan R\,,\qquad\qquad z\doteq 2\arctan S\,,
$$
so that
\beq
R=\tan{w\over 2}\,,\qquad S=\tan{z\over 2}\,.\label{2.9}
\eeq
Using (\ref{2.3}), we obtain the equations
\beq
w_t-c\,w_x={2\over 1+R^2}(R_t-c\,R_x)
={c'\over 2c}{R^2-S^2\over 1+R^2}\,, \label{2.10}
\eeq

\beq
z_t+c\,z_x={2\over 1+S^2}(S_t+c\,S_x)
={c'\over 2c}{S^2-R^2\over 1+S^2}\,. \label{2.11}
\eeq

\begin{figure}[htb]
\centering
\includegraphics[width=0.55\textwidth]{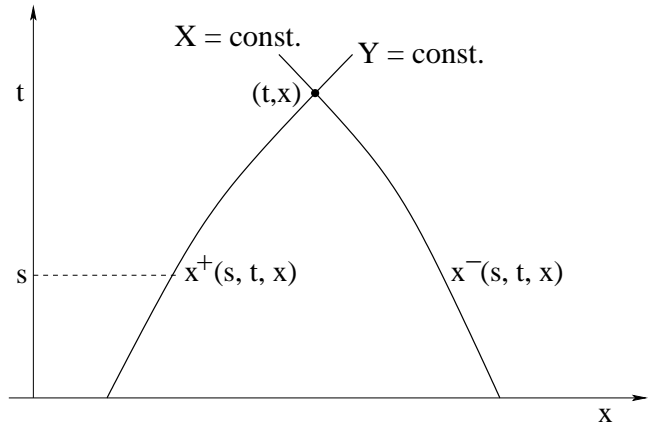} 
\caption{\footnotesize Characteristic curves}
\label{fig1}
\end{figure}

To reduce the equation to a semilinear one, it is convenient to
perform a further change of independent variables (fig.~1).
Consider
the equations for the forward and backward characteristics:
\beq
\dot x^+=c(u)\,,\qquad\qquad \dot x^-=-c(u)\,.\label{2.12} 
\eeq
The characteristics passing through the point $(t,x)$
will be denoted by
$$
s\mapsto x^+(s,t,x)\,,\qquad\qquad s\mapsto x^-(s,t,x)\,,
$$
respectively.
As coordinates $(X,Y)$ of a point $(t,x)$ we shall use the quantities
\beq
X\doteq \int_0^{x^-(0,t,x)} \big(1+R^2(0,x)\big)\,dx
\,,\qquad\qquad Y\doteq \int_{x^+(0,t,x)}^0
\big(1+S^2(0,x)\big)\,dx\,.\label{XY}
\eeq
Of course this implies
\beq
X_t- c(u)X_x =0\,,\qquad\qquad Y_t+ c(u)Y_x =0\,,\label{2.14} 
\eeq
\beq
(X_x)_t- (c\,X_x)_x=0\,,\qquad\qquad (Y_x)_t+(c\,Y_x)_x=0\,.\label{2.15}
\eeq
Notice that
$$
X_x(t,x)= \lim_{h\to 0} {1\over h} \int_{x^-(0, t, x)}
^{x^-(0, t, x+h)}\big(1+ R^2(0, x)\big)\,dx\,,
$$
$$
Y_x(t,x)= \lim_{h\to 0} {1\over h} \int_{x^+(0, t, x)}
^{x^+(0, t, x+h)} \big(1+S^2(0, x)\big)\,dx\,.
$$
For any smooth function $f$, using (\ref{2.14}) one finds
\beq
\begin{array}{ccccccr}
f_t+cf_x &=& f_XX_t+f_Y Y_t+cf_X X_x+cf_Y Y_x & = & (X_t+cX_x)f_X
& = & 2cX_x f_X\,,\\
f_t-cf_x &=& f_XX_t+f_Y Y_t-cf_X X_x-cf_Y Y_x & = & (Y_t-cY_x)f_Y
& = & -2cY_x f_Y\,.
\end{array} \label{2.16} 
\eeq
We now introduce the further variables
\beq
p\doteq {1+R^2\over X_x}\,,\qquad\qquad q\doteq {1+S^2\over -Y_x}\,.
\label{2.17}
\eeq
Notice that the above definitions imply
\beq
(X_x)^{-1}= {p\over 1+R^2}=p\,\cos^2{w\over 2}\,,\qquad\qquad
(-Y_x)^{-1}={q\over 1+S^2}=q\,\cos^2{z\over 2}\,.\label{2.18}
\eeq
{}From (\ref{2.10})-(\ref{2.11}), using (\ref{2.16})-(\ref{2.18}), 
we obtain
$$
2c {(1+S^2)\over q}w_Y={c'\over 2c}
{R^2-S^2\over 1+R^2}\,,\qquad\qquad
2c {(1+R^2)\over p}z_X={c'\over 2c}
{S^2-R^2\over 1+S^2}\,.
$$
Therefore
\beq
\left\{
\begin{array}{ccccc}
w_Y&=&{c'\over 4c^2}\,{R^2-S^2\over 1+R^2}\,{q\over 1+S^2}
&=&{c'\over 4c^2}\,\left( \sin^2{w\over 2}\cos^2{z\over 2}-
\sin^2{z\over 2}\cos^2{w\over 2}\right)\,q\,,\\[4mm]
z_X&=&{c'\over 4c^2}\,{S^2-R^2\over 1+S^2}\,{p\over 1+R^2}
&=&{c'\over 4c^2}\left( \sin^2{z\over 2}\cos^2{w\over 2}-
\sin^2{w\over 2}\cos^2{z\over 2}\right)\,p\,.
\end{array}\right.\label{2.19} 
\eeq
Using trigonometric formulas, the above
expressions can be further simplified as
$$
\left\{
\begin{array}{ccc}
w_Y &= &{c'\over 8c^2}\,( \cos z - \cos w) \,q\,,\\[4mm]
z_X &= &{c'\over 8c^2}\,( \cos w - \cos z) \,p\,.
\end{array}\right.
$$

Concerning the quantities $p,q$, we observe that
\beq
c_x~=~c'\, u_x~=~ c'\,{R-S\over 2c}\,. \label{2.20} 
\eeq
Using again (\ref{2.18}) and (\ref{2.15}) we compute
$$
\begin{array}{ccl}
p_t-c\,p_x &= & (X_x)^{-1}\,2R(R_t-cR_x)-(X_x)^{-2}
\big[(X_x)_t-c(X_x)_x\big](1+R^2)\\[4mm]
&=& (X_x)^{-1}\,2R\,{c'\over 4c}(R^2-S^2)-(X_x)^{-2}[c_xX_x](1+R^2)\\[4mm]
&=& {p\over 1+R^2}2R(R^2-S^2){c'\over 4c}-
{p\over 1+R^2} {c'\over 2c} (R-S)(1+R^2)\\[4mm]
&=& {c'\over 2c}\,{p\over 1+R^2} \big[ S(1+R^2)-R(1+S^2)\big]\,.
\end{array}
$$
$$
\begin{array}{ccl}
q_t+c\,q_x &= & (-Y_x)^{-1}\,2S(S_t-cS_x)-(-Y_x)^{-2}
\big[(-Y_x)_t+c(-Y_x)_x\big](1+S^2)\\[4mm]
&=& (-Y_x)^{-1}\,2S\,{c'\over 4c}(S^2-R^2)-(-Y_x)^{-2}[c_x(Y_x)](1+S^2)\\[4mm]
&=& {q\over 1+S^2}2S(S^2-R^2){c'\over 4c}-
{q\over 1+S^2} {c'\over 2c} (S-R)(1+S^2)\\[4mm]
&=& {c'\over 2c}\,{q\over 1+S^2} \big[ R(1+S^2)-S(1+R^2)\big]\,.
\end{array}
$$
In turn, this yields
\beq
\begin{array}{ccl}
p_Y &=& (p_t-c\,p_x){1\over 2c\,(-Y_x)}=
(p_t-c\,p_x){1\over 2c}\,{q\over 1+S^2}\\[4mm]
&=& {c'\over 4c^2}{S(1+R^2)-R(1+S^2)\over (1+R^2)(1+S^2)}\,pq
={c'\over 4c^2}\,\left[ {S\over 1+S^2}-{R\over 1+R^2}\right]\,pq\\[4mm]
&=& {c'\over 4c^2}\,\left[ \tan{z\over 2}\,\cos^2{z\over 2}-
\tan{w\over 2}\,\cos^2{w\over 2}\right]\,pq={c'\over 4c^2}\,
{\sin z-\sin w\over 2}\,pq\,,
\end{array}   \label{2.21} 
\eeq
\beq
\begin{array}{ccl}
q_X &= & (q_t+c\,q_x){1\over 2c\,(X_x)}=
(q_t+c\,q_x){1\over 2c}\,{p\over 1+R^2}\\[4mm]
&=&{c'\over 4c^2}{R(1+S^2)-S(1+R^2)\over (1+S^2)(1+R^2)}\,pq
={c'\over 4c^2}\,\left[ {R\over 1+R^2}-{S\over 1+S^2}\right]\,pq\\[4mm]
&=&{c'\over 4c^2}\,\left[ \tan{w\over 2}\,\cos^2{w\over 2}-
\tan{z\over 2}\,\cos^2{z\over 2}\right]\,pq={c'\over 4c^2}\,
{\sin w-\sin z\over 2}\,pq\,.
\end{array} \label{2.22}
\eeq
Finally, by (\ref{2.16}) we have
\beq
\begin{array}{ccccccc}
u_X&=&(u_t+cu_x)\,{1\over 2c} {p\over 1+R^2}&=& {1\over 2c}\, {R\over 1+R^2}\,p
&=& {1\over 2c}\, \left(\tan{w\over 2}\,\cos^2{w\over 2}\right)\,p\,, \\[4mm]
u_Y& =&(u_t-cu_x)\,{1\over 2c} {q\over 1+S^2}&=&{1\over 2c}\, {S\over 1+S^2}\,q
&=& {1\over 2c}\, \left(\tan{z\over 2}\,\cos^2{z\over 2}\right)\,q\, \,.
\end{array} \label{2.23} 
\eeq

\vs
Starting with the nonlinear equation (\ref{1.1}),
using $X,Y$ as independent variables we thus obtain a semilinear
hyperbolic system with smooth coefficients for the variables $u,w,z,p,q$.
Using some trigonometric identities, the set of equations (\ref{2.19}),
(\ref{2.21})-(\ref{2.22}) and (\ref{2.23}) can be rewritten as
\beq
\left\{
\begin{array}{ccc}
w_Y&=&{c'\over 8c^2}\,( \cos z - \cos w)\,q\,,\\[4mm]
z_X&=&{c'\over 8c^2}\,( \cos w - \cos z)\,p\,,
\end{array}\right. \label{2.24} 
\eeq
\beq
\left\{
\begin{array}{ccc}
p_Y &= & {c'\over 8c^2}\, \big[\sin z-\sin w\big]\,pq\,,\\[4mm]
q_X &= & {c'\over 8c^2}\, \big[\sin w-\sin z\big]\,pq\,,
\end{array}\right.\label{2.25} 
\eeq
\beq
\left\{ 
\begin{array}{ccc}
u_X &= & {\sin w\over 4c} \, p\,,\\[4mm]
 u_Y &= & {\sin z\over 4c} \, q\,.
\end{array}\right. \label{2.26} 
\eeq
\v
\no{\bf Remark 1.}  The function $u$ can be determined by
using either one of the equations in (\ref{2.26}). One
can easily check that the
two equations are compatible, namely
\beq
\begin{array}{rcl}
u_{XY}
&= &-{\sin w\over 4c^2} \,c' u_Y\,p+{\cos w\over 4c}\, w_Y\,p+{\sin w\over
4c}\,p_Y\\[4mm]
&=&  {c'\over 32\,c^3}\,\big[ -2\sin w\sin z+
\cos w\cos z -\cos^2 w-\sin^2 w + \sin w\sin z\big]\,pq\\[4mm]
&=&{c'\over 32\,c^3}\,\big[ \cos(w-z)-1\big]\,pq \\ [4mm] 
&= & u_{Y\!X} \,.
\end{array} \label{compat}
\eeq
\v
\no{\bf Remark 2.} We observe that the new system is
invariant under translation by $2\pi$
in $w$ and $z$.  Actually, it would be more precise to work with
the variables $w^\dagger\doteq e^{iw}$ and $z^\dagger\doteq e^{iz}$.
However, for simplicity we shall use the variables $w,z$,
keeping in mind that they range on the unit circle
$[-\pi,\pi]$ with endpoints identified.
\v
The system (\ref{2.24})-(\ref{2.26}) must now be supplemented
by non-characteristic boundary conditions,
corresponding to (\ref{1.2}). For this purpose, we observe
that $u_0,u_1$ determine the initial values of the functions
$R,S$ at time $t=0$.
The line $t=0$ corresponds
to a curve $\gamma$ in the $(X, Y)$ plane, say
$$
Y = \vp (X), \qquad X\in\R,
$$
where
$Y\doteq \vp(X)$ if and only if
$$
X=\int_0^x \big(1+R^2(0,x)\big)\,dx\,,\quad
Y=- \int_0^x \big(1+S^2(0,x)\big)\,dx\qquad\hbox{for some}~~x\in\R\,.
$$
We can use the variable $x$ as a parameter along the curve $\gamma$.
The assumptions $u_0\in H^1$, $u_1\in \L^2$ imply
$R,S\in \L^2$; to fix the ideas, let
\beq
\E_0\doteq  {1\over 4}\int \big[R^2(0,x)+S^2(0,x)\big]\,dx <\infty
\,.\label{E.0}
\eeq
The two functions
$$
X(x)\doteq \int_0^x \big(1+R^2(0,x)\big)\,dx \,,\qquad\qquad
Y(x)\doteq \int_x^0 \big(1+S^2(0,x)\big)\,dx
$$
are well defined and absolutely continuous.  Clearly,
$X$ is strictly increasing while $Y$
is strictly decreasing. Therefore, the map $X\mapsto \vp(X)$
is continuous and strictly decreasing.  From (\ref{E.0}) it
follows
\beq
\big|X+\vp(X)\big|\leq 4 \E_0\,.\label{xp}
\eeq
As $(t,x)$ ranges over the domain $[0,\infty[\,\times\R$, the
corresponding variables $(X, Y)$ range over the set
\beq
\Omega^+\doteq \big\{ (X,Y)\,; ~~Y\geq \vp(X)\big\}\,.\label{2.27}
\eeq
Along the curve
$$
\gamma\doteq\big\{ (X,Y)\,;~~Y=\vp(X)\big\}\subset\R^2
$$
parametrized by $x\mapsto \big(X(x), \,Y(x)\big)$,
we can thus assign the boundary data
$(\bar w, \bar z, \bar p, \bar q, \bar u)$$\in \L^\infty$
defined by
\beq
\left\{
\begin{array}{rcl}
{\bar w } &= & 2\arctan R(0,x)\,,\\ {\bar z} &= & 2\arctan S(0,x)\,,
\end{array}\right.\qquad \qquad
\left\{
\begin{array}{rcl}
{\bar p} &\equiv & 1\,,\\
{\bar q} &\equiv & 1\,,
\end{array}\right. \qquad\qquad  {\bar u } =u_0(x)\,. \label{2.28}
\eeq
We observe that the identity
\beq
F = \tan {{\bar w}\over 2}-\tan {{\bar z}\over 2} -2c({\bar u}){\bar u}_x = 0
\label{2.29}
\eeq
is identically satisfied
along $\gamma$. A similar identity holds for $G$.

\section{Construction of integral solutions} 

\setcounter{equation}{0}

Aim of this section is to prove a global existence theorem for the
system (\ref{2.24})-(\ref{2.26}), describing the nonlinear wave equation
in our transformed variables.
\v
\no{\bf Theorem 4.} {\it Let the assumptions in Theorem 1 hold.
Then the corresponding
problem (\ref{2.24})-(\ref{2.26}) with boundary data (\ref{2.28})
has a unique solution, defined for all $(X,Y)\in\R^2$.}
\v
In the following, we shall construct the solution on the domain $\Omega^+$
where $Y\geq \vp(X)$. On the complementary set
$\Omega^-$ where $Y<\vp(X)$, the solution
can be constructed in an entirely similar way.

Observing that all equations (\ref{2.24})-(\ref{2.26})
have a locally Lipschitz continuous right hand side,
the construction of a local solution as fixed
point of a suitable integral transformation is straightforward.
To make sure that this solution is actually defined on the whole
domain $\Omega^+$, one must establish a priori bounds,
showing that $p,q$ remain bounded on bounded sets.
This is not immediately obvious from the equations (\ref{2.25}),
because the right hand sides have quadratic growth.

The basic estimate can be derived as follows. Assume
\beq
C_0\doteq \sup_{u\in \R} \left|{c'(u)\over 4c^2(u)}\right| <\infty\,.
\label{3.1}
\eeq
{}From (\ref{2.25}) it follows the identity
$$
q_X + p_Y = 0\,.
$$
In turn, this implies that the differential form
$p\,dX-q\,dY$ has zero integral
along every closed curve contained in $\Omega^+$.
In particular, for every $(X,Y)\in\Omega^+$,
consider the closed curve $\Sigma$ (see fig.~2) consisting of:

\begin{description}
\item[-] the vertical segment joining $\big( X,\,\vp(X)\big)$ with $(X,Y)$,

\item[-] the horizontal segment joining $(X,Y)$ with 
$\big( \vp^{-1}(Y),\,Y\big)$

\item[-] the portion of boundary $\gamma=\big\{Y=\vp(X)\big\}$ joining
$\big(\vp^{-1}(Y),\,Y\big)$ with $\big(X,\,\vp(X)\big)$.
\end{description}

\begin{figure}[htb]
\centering
\includegraphics[width=0.55\textwidth]{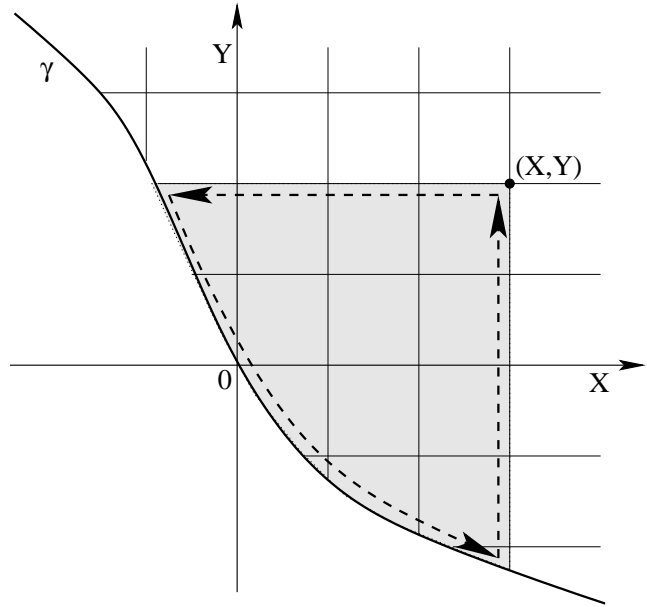} 
\caption{\footnotesize The closed curve $\Sigma$ }
\label{fig2}
\end{figure}

\no Integrating along $\Sigma$, recalling that $p=q=1$ along
$\gamma$ and then using (\ref{xp}), we obtain
\beq
\begin{array}{rcl}
\int_{\vp^{-1}(Y)}^{X}p(X',Y)\,dX' + \int_{\vp(X)}^{Y}q(X, Y')\,dY'
&=& X - \vp^{-1}(Y) + Y -\vp(X) \\[4mm]
&\leq & 2(|X|+|Y|+4\E_0)\,.
\end{array} \label{3.2}
\eeq
Using (\ref{3.1})-(\ref{3.2}) in (\ref{2.25}), since $p,q>0$ we
obtain the a priori bounds
\beq
\begin{array}{rcl}
p(X,Y) &= &\exp\left\{
\int_{\vp(X)}^Y{c'(u)\over 8c^2(u)}\,
\big[\sin z-\sin w\big]\,q(X,Y')\,dY'\right\}\\[4mm]
&\leq&  \exp\left\{ C_0\,
\int_{\vp(X)}^Y\,q(X,Y')\,dY'\right\}\\[4mm]
& \leq & \exp\big\{ 2C_0(|X|+|Y|+ 4\E_0)\big\}\,.
\end{array} \label{3.3}
\eeq
Similarly,
\beq
q(X,Y)\leq \exp \big\{ 2C_0( |X|+|Y|+4\E_0)\big\} \,.\label{3.4}
\eeq

Relying on (\ref{3.3})-(\ref{3.4}), we now show that, on bounded sets
in the $X$-$Y$ plane, the
solution of (\ref{2.24})--(\ref{2.26}) with boundary conditions
(\ref{2.28}) can be obtained as the
fixed point of a contractive transformation.
For any given $r>0$, consider the bounded domain
$$
\Omega_r\doteq \big\{ (X,Y)\,; ~~Y\geq \vp(X)\,,\quad X\leq r\,,
\quad Y\leq r\big\}\,.
$$
Introduce the space of functions
$$
\Lambda_r\doteq\Big\{ f:\Omega_r\mapsto\R\,;~~\|f\|_*\doteq {\rm ess}\!\!
\sup_{(X,Y)\in\Omega_r}~
e^{-\kappa(X+Y)} \big|f(X,Y)\big|~<~\infty\Big\}
$$
where $\kappa$ is a suitably large constant, to be determined later.
For $w,z,p,q,u\in\Lambda_r$, consider the transformation
$\T(w,z,p,q,u)=(\tilde w,\tilde z,\tilde p,\tilde q,\tilde u)$
defined by
\beq
\left\{
\begin{array}{rcl}
\tilde w(X,Y)&=&\bar w(X,\vp(X))+\int_{\vp(X)}^Y
{c'(u)\over 8c^2(u)}\,( \cos z - \cos w )\,q\,dY,\\[4mm]
\tilde z(X,Y)&=&\bar z(\vp^{-1}(Y),Y)+  \int_{\vp^{-1}(Y)}^X
{c'(u)\over 8c^2(u)}\,( \cos w - \cos z)\,p\,dX\,,
\end{array}\right. \label{3.5}
\eeq
\beq
\left\{
\begin{array}{rcl}
\tilde p(X,Y) &=& 1+\int_{\vp(X)}^Y{c'(u)\over 8c^2(u)}\,
\big[\sin z-\sin w\big]\,\hat p\hat q\,dY\},\\[4mm]
\tilde q(X,Y)&=&1+\int_{\vp^{-1}(Y)}^X {c'(u)\over 8c^2(u)}\,
\big[\sin w-\sin z\big]\,\hat p\hat q\,dX\}\,,
\end{array}\right.\label{3.6}
\eeq
\beq
\tilde u(X,Y)=\bar u(X,\vp(X))+\int_{\vp(X)}^Y
{\sin z\over 4c} \, q\,dY.\label{3.7}
\eeq
In (\ref{3.6}), the quantities $\hat p,\hat q$ are defined as
\beq
\hat p\doteq \min\big\{ p\,,~2 e^{2C_0(|X|+|Y|+4\E_0)}\big\}\,,\qquad\qquad
\hat q\doteq \min\big\{ q\,,~2 e^{2C_0(|X|+|Y|+4\E_0)}\big\}\,.\label{3.8}
\eeq
Notice that $\hat p=p$, $\hat q=q$ as long as the a priori estimates
(\ref{3.3})-(\ref{3.4}) are satisfied.  Moreover, if in
the equations (\ref{2.24})-- (\ref{2.26}) the variables $p,q$ are replaced
with $\hat p,\hat q$, then the right hand sides
become uniformly Lipschitz continuous on
bounded sets in the $X$-$Y$ plane.
A straightforward computation now shows that the map $\T$
is a strict contraction on the space $\Lambda_r$, provided that
the constant $\kappa$ is chosen sufficiently big
(depending on the function $c$ and on $r$).

Obviously, if $r'>r$, then the solution of (\ref{3.5})--(\ref{3.7})
on $\Omega_{r'}$ also provides the solution
to the same equations on $\Omega_r$, when restricted to this smaller domain.
Letting $r\to\infty$, in the limit we thus obtain a unique
solution $(w,z,p,q,u)$ of (\ref{3.5})--(\ref{3.7}),
defined on the whole domain $\Omega^+$.

To prove that these functions satisfy the  (\ref{2.24})--(\ref{2.26}),
we claim that $\hat p=p$, $\hat q=q$ at every point $(X,Y)\in\Omega^+$.
The proof is by contradiction.
If our claim does not hold, since the maps $Y\mapsto p(X,Y)$,
$X\mapsto q(X,Y)$
are continuous, we can find some point
$(X^*, Y^*)\in\Omega^+$ such that
\beq
p(X,Y)\leq
2 e^{2C_0(|X|+|Y|+4\E_0)}\,,\qquad \qquad q(X,Y)\leq 2 e^{2C_0(|X|+|Y|+4\E_0)}
\label{3.9}
\eeq
for all $(X,Y)\in\Omega^*\doteq \Omega^+\cap \{X\leq X^*\,,~~Y\leq Y^*\}$, but
either
$p(X^*,Y^*)\geq {3\over 2}  e^{2C_0(|X|+|Y|+4\E_0)}$ or
$q(X^*,Y^*)\geq {3\over 2} e^{2C_0(|X|+|Y|+4\E_0)}$.
By (\ref{3.9}), we still have $\hat p=p$, $\hat q=q$ restricted to $\Omega^*$,
hence the equations (\ref{2.24})--(\ref{2.26})  and the a priori bounds
(\ref{3.3})-(\ref{3.4}) remain valid.  In particular, these imply
$$
p(X^*,Y^*)\leq e^{2C_0(|X|+|Y|+4\E_0)}\,,\qquad\qquad q(X^*,Y^*)
\leq e^{2C_0(|X|+|Y|+4\E_0)}\,,
$$
reaching a contradiction.
\endproof

\v
\no{\bf Remark 3.} In the solution constructed above,
the variables $w, z$ may well grow outside the initial range $\,]-\pi, \pi[\,$.
This happens precisely when the quantities $R$, $S$
become unbounded, i.e.~when
singularities arise.
\v
For future reference, we state a useful
consequence of the above construction.
\v
\no{\bf Corollary 1.} {\it If the initial data $u_0, u_1$ are smooth,
then the solution $(u,p,q,w,z)$ of (\ref{2.24})--(\ref{2.26}), (\ref{2.28})
is a smooth function of the variables $(X,Y)$.
Moreover, assume that a sequence of smooth functions
$(u_0^m, u_1^m)_{m\geq 1}$ satisfies
$$
u_0^m\to u_0\,,\qquad (u_0^m)_x\to (u_0)_x\,,\qquad u_1^m\to u_1
$$
uniformly on compact subsets of $\R$. Then one
has the convergence of the corresponding
solutions:
$$(u^m,\,p^m,\, q^m,\, w^m,\, z^m)\to (u,p,q,w,z)$$
uniformly on bounded subsets of the $X$-$Y$ plane.}
\v
We also remark that the equations (\ref{2.24})--(\ref{2.26})
imply the conservation laws
\beq
q_X+p_Y=0\,,\qquad\qquad \left({q\over c}\right)_X
-\left({p\over c}\right)_Y=0\,.\label{3.10}
\eeq

\section{Weak solutions, in the original variables} 

\setcounter{equation}{0}

By expressing the solution $u(X,Y)$ in terms of the original variables
$(t,x)$, we shall recover a solution of the Cauchy problem
(\ref{1.1})-(\ref{1.2}). This will provide a proof of Theorem 1.

As a preliminary, we examine the regularity of the solution
$(u, w, z, p, q)$ constructed in the previous section.
Since the initial data $(u_0)_x$ and $u_1$ are only assumed to
be in $\L^2$, the functions $w,z,p,q$ may well be discontinuous.
More precisely, on bounded subsets of the $X$-$Y$ plane,
the equations (\ref{2.24})-(\ref{2.26}) imply
the following:
\v
\begin{description}

\item[-] The functions $w,p$ are Lipschitz continuous w.r.t.~$Y$,
measurable w.r.t.~$X$.

\item[-] The functions $z,q$ are Lipschitz continuous w.r.t.~$X$,
measurable w.r.t.~$Y$.

\item[-] The function $u$ is Lipschitz continuous w.r.t.~both $X$ and $Y$.
\end{description}
\v
The map $(X,Y)\mapsto (t,x)$ can be constructed as follows.
Setting $f=x$, then $f=t$  in the two equations at
(\ref{2.16}), we find
$$
\left\{
\begin{array}{rcr}
c &=& 2cX_x\,x_X\,, \\[3mm]
-c&=& -2cY_x\,x_Y\,,
\end{array}\right.\qquad\qquad
\left\{
\begin{array}{rcr}
1&= &2cX_x\,t_X\,,\\[3mm]
1&=& -2cY_x\,t_Y\,,
\end{array}\right.
$$
respectively.
Therefore, using (\ref{2.18}) we obtain
\beq
\left\{
\begin{array}{ccccr} 
x_X&~=~&{1\over 2X_x}&~=~&{(1+\cos w)\,p\over 4}\,,\\[4mm]
x_Y&~=~&{1\over 2Y_x}&~=~&-{(1+\cos z)\,q\over 4}\,,
\end{array}\right. \label{4.1}
\eeq
\beq
\left\{
\begin{array}{ccccr} 
 t_X&=&{1\over 2cX_x}&=&{(1+\cos w)\,p\over 4c}\,,\\[4mm]
t_Y&=&{1\over -2cY_x}&=&{(1+\cos z)\,q\over 4c}\,.
\end{array}\right. \label{4.2}
\eeq
For future reference,
we write here the partial derivatives of the inverse
mapping, valid at points where $w,z\not= -\pi$.
\beq
\left\{
\begin{array}{ccr}  
X_x&=&{2\over (1+\cos w)\,p}\,,\\[4mm]
Y_x&=&\,-{2\over(1+\cos z)\,q}\,,
\end{array}\right. 
\qquad\qquad
\left\{
\begin{array}{ccr} 
 X_t&=&{2c\over (1+\cos w)\,p}\,,\\[4mm]
Y_t&=&{2c\over(1+\cos z)\,q}\,.
\end{array}\right. \label{4.3}
\eeq
We can now recover the functions $x=x(X,Y)$ by integrating
one of the  equations in (\ref{4.1}). Moreover, we can
compute $t=t(X,Y)$ by integrating one of the equations in (\ref{4.2}).
A straightforward calculation shows that the two equations
in (\ref{4.1}) are equivalent: differentiating the first w.r.t.~$Y$
or the second w.r.t.~$X$ one obtains the same expression.
$$
\begin{array}{ccl} 
 x_{XY} & = &{(1+\cos w)\,
p_Y\over 4} -{p\sin w\, w_Y\over 4}\\[3mm]
  & = & {c'\,pq\over 32 c^2}\,\big[ \sin z-\sin w+\sin(z-w)\big]
=x_{Y\!X}\,.
\end{array}
$$
Similarly, the equivalence of the two equations in (\ref{4.2})
is checked by
$$
\begin{array}{cl} 
&t_{XY}-t_{Y\!X}~=~
\left({x_X\over c}\right)_Y+\left({x_Y\over c}\right)_X~=~{2\over c}\, x_{XY}
-\left({x_X\over c^2}\,c' u_Y+{x_Y\over c^2}\,c'u_X\right)\\[3mm]
&\quad={c'\,pq\over 16\, c^3}\,\big[ \sin z-\sin w+\sin(z-w)\big]\\[3mm]
&\qquad\qquad
 -{c'\,pq\over 16\, c^3}\big[ (1+\cos w)\sin z-(1+\cos z)\sin w\big]~=~0\,.
\end{array}
$$
\v
In order to define $u$ as a function of the original variables $t,x$,
we should formally invert the map $(X,Y)\mapsto (t,x)$ and write
$u(t,x)=u\big( X(t,x)\,,~Y(t,x)\big)$.
The fact that
the above map may  not be one-to-one does not
cause any real difficulty. Indeed, given $(t^*,x^*)$, we can choose an
arbitrary point $(X^*,Y^*)$ such that $t(X^*,Y^*)=t^*$, $x(X^*,Y^*)=x^*$,
and define $u(t^*,x^*)=u(X^*,Y^*)$.
To prove that the values of $u$ do not depend on the choice
of $(X^*,Y^*)$, we
proceed as follows.
Assume that there are two distinct points such that
$t(X_1, Y_1)=t(X_2,Y_2)=t^*$,
$~x(X_1,Y_1)=x(X_2,Y_2)=x^*$.
We consider two cases:
\v
\no Case 1: $X_1\leq X_2$, $Y_1\leq Y_2$. Consider the set
$$
\Gamma_{x^*}\doteq \Big\{ (X,Y)\,;~~x(X,Y)\leq x^*\Big\}
$$
and call $\partial \Gamma_{x^*}$ its boundary.
By (\ref{4.1}), $x$ is increasing with $X$ and decreasing with $Y$.
Hence, this boundary can be represented as
the graph of a Lipschitz continuous function: $X-Y=\phi(X+Y)$.
We now construct the Lipschitz continuous curve $\gamma$ (fig.~3a) consisting
of
\begin{description}

\item[-] a horizontal segment joining $(X_1,Y_1)$ with a point $A=(X_A,Y_A)$
on $\partial \Gamma_{x^*}$, with $Y_A=Y_1$,

\item[-] a portion of the boundary $\partial \Gamma_{x^*}$,

\item[-] a vertical segment joining $(X_2,Y_2)$ to a point $B=(X_B,Y_B)$
on $\partial \Gamma_{x^*}$, with $X_B=X_2$.
\end{description}
\no We can obtain a Lipschitz continuous parametrization
of the curve $\gamma:[\xi_1,\xi_2]\mapsto \R^2$ in terms of
the parameter $\xi=X+Y$.
Observe that the map $(X,Y)\mapsto (t,x)$ is constant along $\gamma$.
By (\ref{4.1})-(\ref{4.2}) this implies
$(1+\cos w) X_\xi=(1+\cos z)Y_\xi = 0$,
hence $\sin w\cdot X_\xi=\sin z\cdot Y_\xi=0$. We now compute
$$
\begin{array}{rl}
u(X_2,Y_2)-&u(X_1,Y_1)=\int_\gamma \big(u_X\,dX+u_Y\,dY\big) \cr
&=\int_{\xi_1}^{\xi_2} \left({p\,\sin w\over 4c}\,X_\xi-{q\,\sin z\over 4c}\,
Y_\xi\right)\,d\xi=0\,,
\end{array}
$$
proving our claim.
\v
\no Case 2: $X_1\leq X_2$, $Y_1\geq Y_2$.   In this case, we consider the
set
$$
\Gamma_{t^*}\doteq \Big\{ (X,Y)\,;~~t(X,Y)\leq t^*\Big\}\,,
$$
and construct a curve $\gamma$ connecting $(X_1,Y_1)$ with
$(X_2,Y_2)$ as in fig.~3b.  Details are entirely similar to Case 1.

\begin{figure}[htb]
\centering
\includegraphics[width=0.55\textwidth]{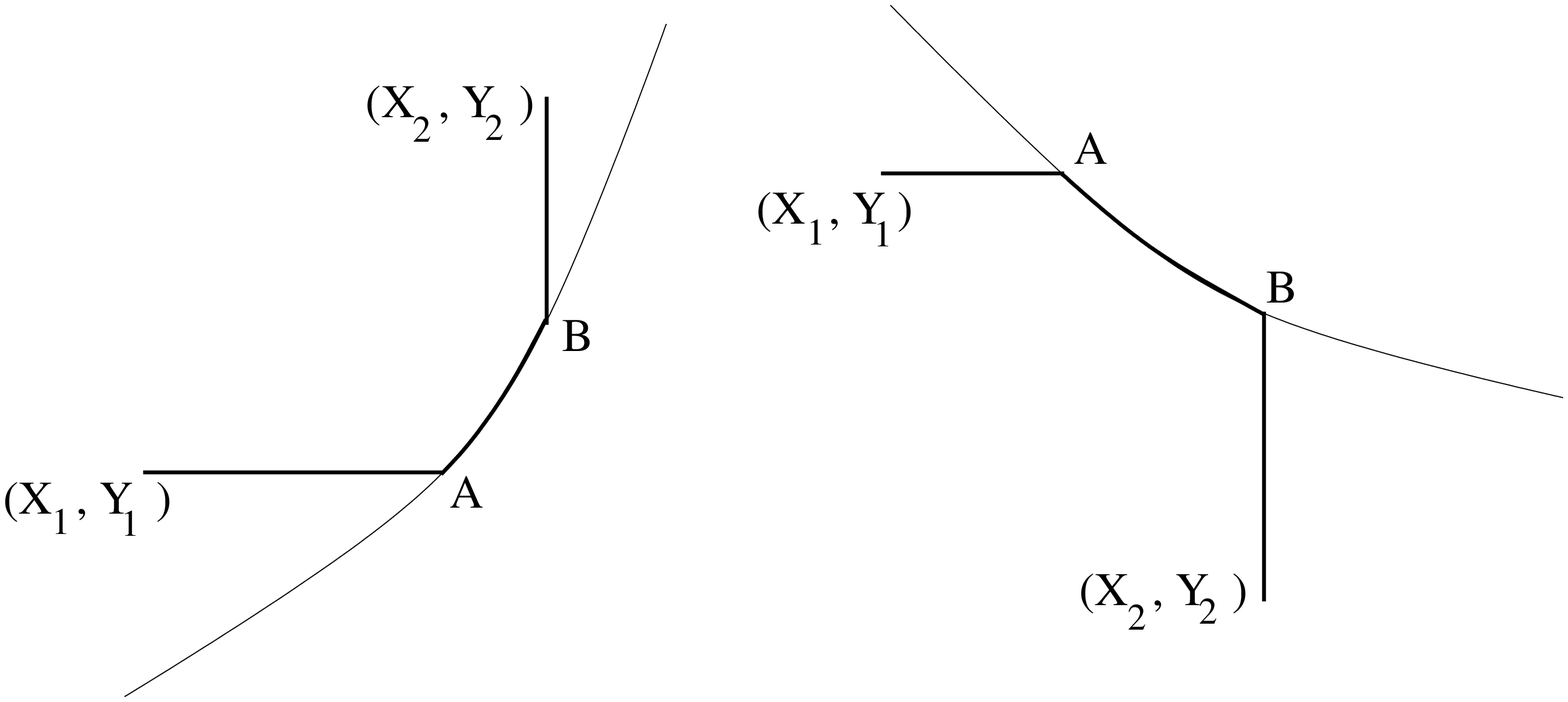} 
\caption{\footnotesize Paths of integration }
\label{fig3}
\end{figure}

We now prove that the function $u(t,x)=u\big(X(t, x), Y(t, x)\big)$
thus obtained is H\"older continuous on bounded sets.
Toward this goal, consider any characteristic curve, say  $t\mapsto
x^+(t)$, with $\dot x^+=c(u)$.   By construction, this
is parametrized by the function $X\mapsto \big(t(X,\ov Y),\, x(X,\ov Y)\big)$,
for some fixed $\ov Y$. Recalling (\ref{2.16}), (\ref{2.14}), (\ref{2.18}) 
and (\ref{2.26}), we  compute
\beq
\begin{array}{rcl}
 \int_0^\tau \big[ u_t+c(u)u_x\big]^2\,dt &= &\int_{X_0}
^{X_\tau} ( 2cX_x u_X)^2\, (2X_t)^{-1} dX\\[4mm]
&=& \int_{X_0}
^{X_\tau}  2c\left( p\,\cos^2{w\over 2}\right)^{-1}\,
\left({p\over 4c}\,2\sin{w\over 2}\cos {w\over 2}\right)^2\,dX \\[4mm]
&\leq& \int_{X_0}
^{X_\tau}{p\over 2c}\,dX ~\leq ~C_\tau\,,
\end{array}\label{4.4}
\eeq
for some constant $C_\tau$ depending only on $\tau$.
Similarly, integrating along any backward characteristics $t\mapsto x^-(t)$
we obtain
\beq
\int_0^\tau \big[ u_t-c(u)u_x\big]^2\,dt\leq C_\tau. \label{4.5}
\eeq
Since the speed of characteristics is $\pm c(u)$,
and $c(u)$ is uniformly positive and
bounded, the bounds (\ref{4.4})-(\ref{4.5})
imply that the function $u=u(t,x)$ is H\"older
continuous with exponent $1/2$.  In turn, this implies that
all characteristic curves are $\C^1$ with H\"older continuous derivative.
Still from (\ref{4.4})-(\ref{4.5}) it follows that the functions
$R,S$ at (\ref{2.1}) are square integrable on bounded subsets of the
$t$-$x$ plane.
\v
Finally, we prove that the function $u$ provides a weak solution
to the nonlinear wave equation (\ref{1.1}). According to (\ref{1.4'}),
we need to show that
\beq
\begin{array}{rcl}
0&= &\dint \phi_t\big[(u_t+cu_x)+(u_t-cu_x)\big]
-\big(c(u)\phi\big)_x\big[(u_t+cu_x)-(u_t-cu_x)\big]
\,dxdt\\[4mm]
&=&\dint \big[\phi_t
-(c\phi)_x\big]\,(u_t+cu_x)\,dxdt+\dint \big[\phi_t
+(c\phi)_x\big]\,(u_t-cu_x)\,dxdt\\[4mm]
&=&\dint \big[\phi_t
-(c\phi)_x\big]\,R\,dxdt+\dint \big[\phi_t
+(c\phi)_x\big]\,S\,dxdt\,.
\end{array}\label{4.6}
\eeq
By (\ref{2.16}), this is equivalent to
\beq
\dint \Big\{-2cY_x \phi_Y \,R~+~2cX_x\phi_X \,S~+~
c'(u_X\,X_x+u_Y\,Y_x) \,\phi\,(S-R)\Big\} \,dxdt =0\,.
\label{4.7}
\eeq
It will be convenient to express the double integral in (\ref{4.7})
in terms of the variables $X,Y$. We notice that, by (\ref{2.18}) and 
(\ref{2.14}),
$$
dx\,dt={pq\over 2c\, (1+R^2)(1+S^2)}\,dXdY\,.
$$
Using (\ref{2.26}) and the identities
\beq
\left\{
\begin{array}{ccc}
{1\over 1+R^2}&=&\cos^2{w\over 2}={1+\cos w\over 2}\,, \\[4mm]
{1\over 1+S^2}&=&\cos^2{z\over 2}={1+\cos z\over 2}\,,
\end{array}\right.
\qquad\qquad
\left\{
\begin{array}{ccc}
{R\over 1+R^2} &=&{\sin w\over 2}\,,\\[4mm]
{S\over 1+S^2} &=&{\sin z\over 2}\,,
\end{array}\right.\label{4.8}
\eeq
the double integral in (\ref{4.6}) can thus be written as
\beq
\begin{array}{rl}
\dint &\left\{ 2c\,{1+S^2\over q} \,\phi_Y\,R+
2c{1+R^2\over p}\,\phi_X\,S
+c'\left( {\sin w\over 4c}\,p\,{1+R^2\over p}-
{\sin z\over 4c}\,q\,{1+S^2\over q}\right)
\,\phi\,(S-R)\right\}\\[4mm]
&\qquad\qquad\cdot {pq\over 2c\,(1+R^2)\,(1+S^2)}\,dXdY\\[4mm]
&=\dint \left\{ {R\over 1+R^2}\,p\,\phi_Y+{S\over 1+S^2}\,q\,\phi_X
+{c'pq\over 8c^2}\left({\sin w\over 1+S^2}-{\sin z\over 1+R^2}
\right)\phi\,(S-R)\right\}\,dXdY\\[4mm]
&=\dint\bigg\{ {p\,\sin w\over 2}\,\phi_Y+{q\sin z\over 2}\,\phi_X\\[4mm]
&\qquad +
{c'pq\over 8c^2}\left( \sin w\sin z-\sin w\,\cos^2 {z\over 2}\tan{w\over 2}
-\sin z\,
\cos^2 {w\over 2}\,
\tan{z\over 2}\right)\,\phi\bigg\}\,dXdY\\[4mm]
&= \dint\left\{ {p\,\sin w\over 2}\,\phi_Y+{q\sin z\over 2}\,\phi_X +
{c'pq\over 8c^2}\,\big[\cos (w+z)-1\big]\,\phi\right\} \,dXdY\,.
\end{array} \label{4.9}
\eeq
Recalling (\ref{2.27}), one finds
\beq
\begin{array}{rcl}
\left(p\,\sin w\over 2\right)_Y+
\left(q\,\sin z\over 2\right)_X
&=&(2c\,u_X)_Y+(2c\,u_Y)_X\\[4mm]
&=&4c' u_Xu_Y+4c\,u_{XY}\\[4mm]
&=&{c'\,pq\over 4c^2}\, \sin w\,\sin z+ {c'\,pq\over 8c^2}
\,\big[\cos(w-z)-1\big]\\[4mm]
&=& {c'\,pq\over 8c^2}\big[\cos(w+z)-1\big]\,.
\end{array}\label{4.10}
\eeq
Together, (\ref{4.9}) and (\ref{4.10}) imply (\ref{4.7}) and hence (\ref{4.6}).
This establishes the integral equation (\ref{1.4'}) for every
test function $\phi\in\C^1_c$.

\section{Conserved quantities} 

\setcounter{equation}{0}

From the conservation laws (\ref{3.10}) it follows that the 1-forms
$p \,dX-q\,dY$ and $\f{p}{c}\,dX + \f{q}{c}\,dY$ are closed, hence 
their integrals along any closed curve in th $X$-$Y$ plane vanish.
From the conservation laws at (\ref{2.6}), it follows that
the 1-forms
\beq
E\,dx-(c^2M)\,dt\,,\qquad\qquad M\,dx-E\,dt \label{5.1}
\eeq
are also closed. There is a simple correspondence. In fact
$$
E\,dx-(c^2M)\,dt = \f{p}{4}\,dX -\f{q}{4}\,dY -\f12\,dx
\quad 
-M\,dx+ E\,dt  = \f{p}{4c}\,dX + \f{q}{4c}\,dY -\f12\,dt.
$$
Recalling (\ref{4.1})-(\ref{4.2}), these can be written
in terms of the $X$-$Y$ coordinates as
\beq
{(1-\cos w)\,p\over 8}\,dX-{(1-\cos z)\,q\over 8}\,dY\,,\label{5.2}
\eeq
\beq
{(1-\cos w)\,p\over 8c}\,dX+{(1-\cos z)\,q\over 8c}\,dY\,,\label{5.3}
\eeq
respectively.
Using (\ref{2.24})--(\ref{2.26}), one easily checks that these forms are indeed
closed:
\beq
\left({(1-\cos w)\,p\over 8}\right)_Y={c'pq\over 64c^2}\Big[
\sin z(1-\cos w)-\sin w(1-\cos z)\Big]
=-\left({(1-\cos z)\,q\over 8}\right)_X\,,\label{sour}
\eeq
$$
\left({(1-\cos w)\,p\over 8c}\right)_Y=
{c'pq\over 64c^3}\Big[\sin(w+z)-(\sin w +\sin z)\Big]
=\left({(1-\cos z)\,q\over 8c}\right)_X\,.
$$
In addition, we have the 1-forms
\beq
dx={(1+\cos w)\,p\over 4}\,dX-{(1+\cos z)\,q\over 4}\,dY\,,\label{5.4}
\eeq
\beq
dt={(1+\cos w)\,p\over 4c}\,dX+{(1+\cos z)\,q\over 4c}\,dY\,,\label{5.5}
\eeq
which are obviously closed.

\begin{figure}[htb]
\centering
\includegraphics[width=0.55\textwidth]{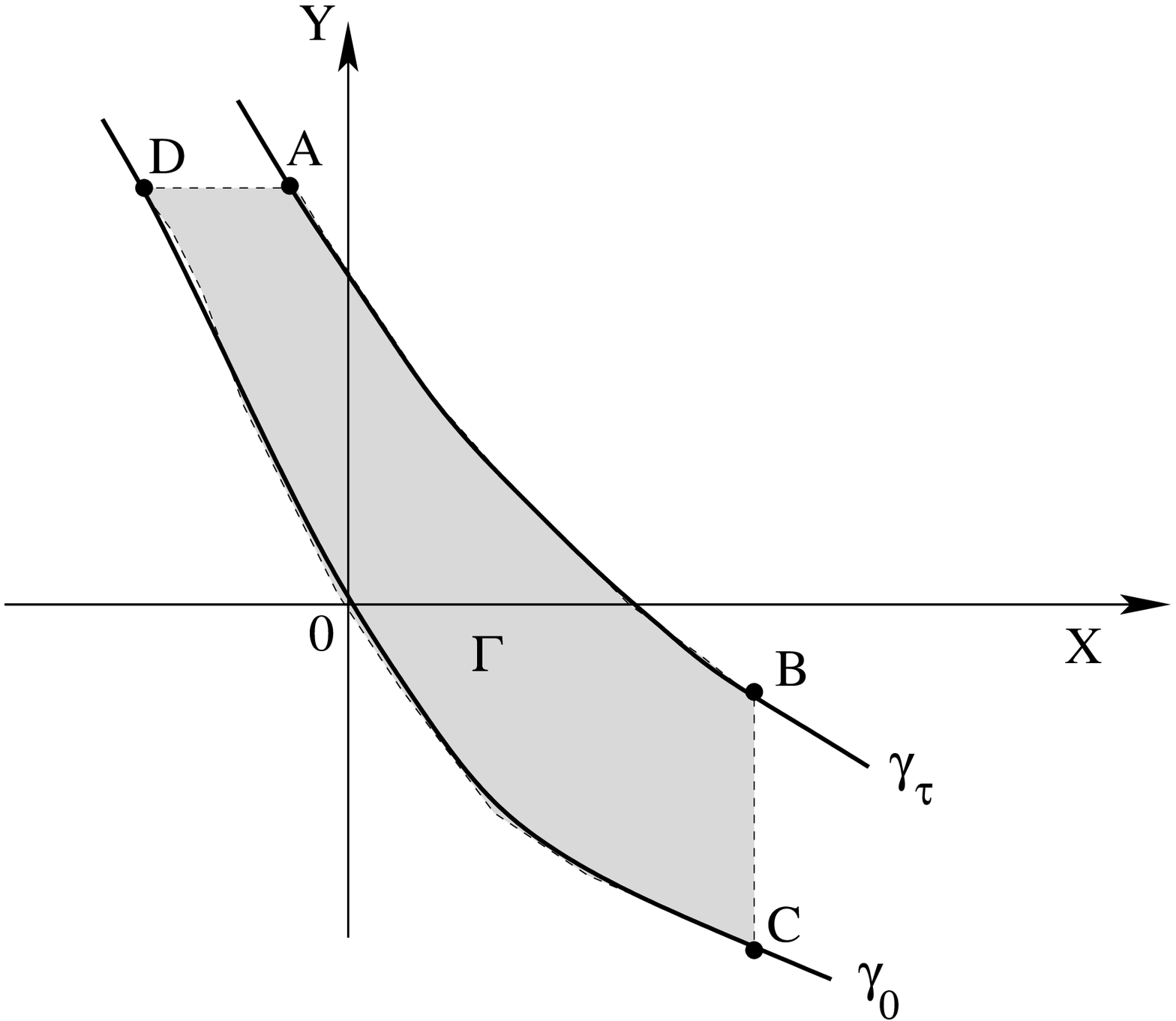} 
\caption{\footnotesize  }
\label{fig4}
\end{figure}

The solutions $u=u(X,Y)$ constructed in Section 3 are
{\it conservative}, in the sense that the integral of
the form (\ref{5.2}) along every Lipschitz continuous, closed curve
in the $X$-$Y$ plane is zero.     

To prove the inequality (\ref{1.7}), fix any $\tau>0$.  
The case $\tau<0$ is identical. For a given $r>0$ arbitrarily large,
define the set (fig.~4)
\beq
\Gamma\doteq \Big\{ (X,Y)\,;~~0\leq t(X,Y)\leq \tau\,,\quad X\leq r\,,
~~Y\leq r\Big\}\label{5.6}
\eeq
By construction, the the map $(X,Y)\mapsto (t,x)$ will act as follows:
$$
A\mapsto (\tau, a)\,,\qquad B\mapsto (\tau, b)\,,\qquad C\mapsto(0,c)\,,
\qquad D\mapsto (0,d)\,,
$$
for some $a<b$ and $d<c$.
Integrating the 1-form (\ref{5.2}) along the boundary of $\Gamma$ we
obtain
\beq
\begin{array}{rl}
\int_{AB}& {(1-\cos w)\,p\over 8}\,dX-{(1-\cos z)\,q\over 8}\,dY\, \\[4mm]
&=\int_{DC} {(1-\cos w)\,p\over 8}\,dX-{(1-\cos z)\,q\over 8}\,dY\,
-\int_{DA}{(1-\cos w)\,p\over 8}\,dX-\int_{CB}{(1-\cos z)\,q\over 8}\,dY\\[4mm]
&\leq \int_{DC} {(1-\cos w)\,p\over 8}\,dX-{(1-\cos z)\,q\over 8}\,dY\,
\\[4mm]
&=\int_d^c \f12\Big[ u_t^2(0,x)+c^2\big(u(0,x)\big)\,u_x^2(0,x)\Big]\,dx
\,.
\end{array}
\label{5.7}
\eeq
On the other hand, using (\ref{5.4}) we compute
\beq
\begin{array}{rl}
\int_a^b \f12 \Big[ u_t^2(\tau,x)+& c^2\big(u(\tau,x)\big)\,u_x^2(\tau,x)
\big]\,dx\\[4mm]
&=\int_{AB\cap\{\cos w\not= -1\}}
{(1-\cos w)\,p\over 8}\,dX+\int_{AB\cap \{\cos z\not= -1\}}
{(1-\cos z)\,q\over 8}\,dY\\[4mm]
&\leq \E_0\,.
\end{array}\label{5.8}
\eeq
Notice that the last relation in (\ref{5.7}) is satisfied as an equality,
because at time $t=0$, along the curve $\gamma_0$
the variables  $w,z$ never assume the value $-\pi$.
Letting $r\to +\infty$ in (\ref{5.6}), one has
$a\to -\infty$, $b\to +\infty$.
Therefore (\ref{5.7}) and (\ref{5.8}) together imply
$\E(t)\leq \E_0$, proving (\ref{1.7}).
\v

\begin{figure}[htb]
\centering
\includegraphics[width=0.95\textwidth]{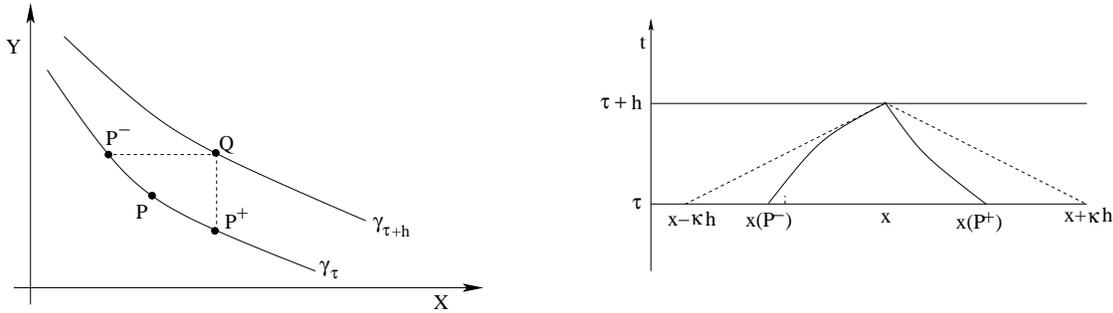} 
\caption{\footnotesize Proving Lipschitz continuity }
\label{fig5}
\end{figure}

We now prove the Lipschitz continuity of the map $t\mapsto
u(t,\cdot)$ in the $\L^2$ distance.
For this purpose, for any fixed time $\tau$, we let
$\mu_\tau=\mu_\tau^- + \mu_\tau^+$
be the positive measure on the real line
defined as follows. 
In the smooth case,
\beq
\mu_\tau^-\big(\,]a,b[\,\big)
=\frac{1}{4}\int_a^b R^2(\tau,x)\,dx\,,\qquad\qquad
\mu_\tau^+\big(\,]a,b[\,\big)
=\frac{1}{4}\int_a^b S^2(\tau,x)\,dx\,.
\label{5.9'}
\eeq
To define $\mu_\tau^\pm$ in the general case,
let $\gamma_\tau$ be the boundary of the
set
\beq
\Gamma_\tau\doteq \big\{(X,Y)\,;~~t(X,Y)\leq\tau\big\}\,.\label{5.9}
\eeq
Given any open interval $]a,b[\,$,
let $A=(X_A, Y_A)$ and $B=(X_B, Y_B)$ be the points on
$\gamma_\tau$ such that
$$
x(A)=a\,,\qquad X_P-Y_P\leq X_A-Y_A\quad\hbox{for every
point $P\in \gamma_\tau$ with $x(P)\leq a$}\,,
$$
$$
x(B)=b\,,\qquad X_P-Y_P\geq X_B-Y_B\quad\hbox{for every
point $P\in \gamma_\tau$ with $x(P)\geq b$}\,.
$$
Then
\beq
\mu_\tau\big(\,]a,b[\,\big)
=\mu_\tau^-\big(\,]a,b[\,\big)+\mu_\tau^+\big(\,]a,b[\,\big)\,,
\label{5.10}
\eeq
where
\beq
\mu_\tau^-\big(\,]a,b[\,\big)\doteq \int_{AB}
{(1-\cos w)\,p\over 8}\,dX\,\qquad
\mu_\tau^+\big(\,]a,b[\,\big)\doteq -\int_{AB}
{(1-\cos z)\,q\over 8}\,dY\,.
\label{5.10'}
\eeq

Recalling the discussion at (\ref{5.1})--(\ref{5.2}), it is clear 
that $\mu_\tau^-$, $\mu_\tau^+$
are bounded, positive measures, and
$\mu_\tau(\R)=\E_0$, for all $\tau$.
Moreover, by (\ref{5.9'}) and (\ref{2.5}),
$$
\int_a^b c^2 u_x^2\,dx \leq \int_a^b\f12(R^2+S^2)dx \leq 2\mu(]a, b[).
$$
For any $a < b$, this yields the estimate
\beq
|u(\tau, b) -u(\tau, a)|^2 \leq |b-a|\int_a^b u_x^2(\tau, y)\,dy\leq 
2\kappa^2|b-a|\mu_\tau(]a, b[). 
\label{5.11}
\eeq
Next, for a given $ h>0$, $y\in\R$, we seek an estimate on the distance
$\big| u(\tau+ h, \,y)-u(\tau,y)|$. As in fig.~5,
let $\gamma_{\tau+ h}$ be the boundary of the set
$\Gamma_{\tau+ h}$, as in (\ref{5.9}).
Let $P=(P_X,\,P_Y)$ be the point on $\gamma_\tau$
such that
$$
x(P)=y\,,\qquad X_{P'}-Y_{P'}\leq X_P-Y_P\quad\hbox{for every
point $P'\in \gamma_\tau$ with $x(P')\leq x$}\,,
$$
Similarly, let $Q=(Q_X,\,Q_Y)$ be the point on $\gamma_{\tau+ h}$
such that
$$
x(Q)=y\,,\qquad X_{Q'}-Y_{Q'}\leq X_Q-Y_Q\quad\hbox{for every
point $Q'\in \gamma_{\tau+ h}$ with $x(Q')\leq y$}\,.
$$
Notice that $X_P\leq X_Q$ and $Y_P\leq Y_Q$.
Let $P^+=(X^+, Y^+)$ be a point on $\gamma_\tau$ with
$X^+ = X_Q$, and let $P^-=(X^-, Y^-)$ be a point on $\gamma_\tau$ with
$Y^- = Y_Q$.
Notice that $x(P^+)\in \,]y,\, y+\kappa h[\,$, because the
point $\big(\tau, x(Q)\big)$ lies on some characteristic
curve with speed $-c(u)>-\kappa$, passing through the point $(\tau+ h,\,y)$.
Similarly, $x(P^-)\in \,]y-\kappa h,\,y[\,$. Recalling that the forms 
in (\ref{5.2}) and (\ref{5.5}) are closed, we obtain the estimate
\beq
\begin{array}{rcl}
\big|u(Q)-u(P^+)\big|&\leq&\int_{Y^+}^{Y_Q}
\big|u_Y(X_Q,Y)\big|\,dY\\[4mm]
&=&\int_{Y^+}^{Y_Q}
\left|{\sin z\over 4c}\,q\right|\,dY\\[4mm]
&=&\int_{Y^+}^{Y_Q}
\left({(1+\cos z)\,q\over 4c}\right)^{1/2}\left({(1-\cos z)\,q
\over 4}\right)^{1/2}\,dY\\[4mm]
&\leq &\left(\int_{Y^+}^{Y_Q}
{(1+\cos z)\,q\over 4c}\,dY\right)^{1/2}\cdot
\left(\int_{Y^+}^{Y_Q}{(1-\cos z)\,q
\over 4}\,dY\right)^{1/2}\\[4mm]
&\leq& \left\{ \int_{P^-P^+} \left[
{(1-\cos w)\,p\over 4}\,dX-{(1-\cos z)\,q\over 4}\,dY\right]\right\}^{1/2}
\cdot  h^{1/2}\,.
\end{array}\label{5.12}
\eeq
The last term in (\ref{5.12}) contains the integral of the 1-form at 
(\ref{5.2}),
along the curve $\gamma_\tau$, between $P^-$ and $P^+$. Recalling the
definition (\ref{5.10})--(\ref{5.10'}) and the estimate (\ref{5.11}),
 we obtain the bound
\beq
\begin{array}{rcl}
\big|u(\tau+ h , x)-u(\tau ,x)\big|^2&\leq &
2\big|u(Q)-u(P^+)\big|^2+2\big|u(P^+)-u(P)\big|^2\\[3mm]
&\leq & 4 h\cdot\mu_\tau\big(]x-\kappa h,\,x+\kappa h[\big)
+ 4\kappa^2\cdot(\kappa h)\cdot\mu_\tau \big(]x\,,~x+ h[\big)\,.
\end{array}\label{5.13}
\eeq
Therefore, for any $ h>0$,
\beq
\begin{array}{rcl}
\big\|u(\tau+ h,\cdot)-u(\tau,\cdot)\big\|_{\L^2}
&=&\left(\int\big|u(\tau+ h,x)-u(\tau,x)\big|^2\,dx\right)^{1/2}\\[4mm]
&\leq&\left(\int
4(1+\kappa^3) h \cdot \mu_\tau\big(]x-\kappa h\,,~x+\kappa h[\big)
\,dx\right)^{1/2}\\[4mm]
&= &\Big( 4(\kappa^3+1) h^2\,\mu_\tau(\R)\Big)^{1/2}\\[4mm]
&= & h\cdot\big[4(\kappa^3+1)\,\E_0\big]^{1/2}.
\end{array} \label{5.14}
\eeq
This proves the uniform Lipschitz continuity of the map
$t\mapsto u(t,\cdot)$, stated at (\ref{1.lip}).

\section{Regularity of trajectories}  

\setcounter{equation}{0}

In this section we prove the continuity of the
functions $t\mapsto u_t(t,\cdot)$ and $t\mapsto u_x(t,\cdot)$,
as functions with values in $\L^p$.  This will complete the proof
of Theorem 1.

We first consider the case where the initial data $(u_0)_x$, $u_1$
are smooth with compact support.
In this case, the solution $u=u(X,Y)$ remains smooth on the entire
$X$-$Y$ plane.   Fix a time $\tau$ and let $\gamma_\tau$ be the
boundary of the set $\Gamma_\tau$, as in (\ref{5.9}).
We claim that
\beq
{d\over dt} u(t,\cdot)\bigg|_{t=\tau}= u_t(\tau,\cdot)\label{6.1}
\eeq
where, by (\ref{2.14}), (\ref{2.18}) and (\ref{2.26}),
\beq
\begin{array}{rcl}
u_t (\tau,x) & \doteq&  u_X\,X_t+u_Y\,Y_t \\ [4mm]
 & =&  {p\,\sin w\over 4c}\,
{2c\over p(1+\cos w)}+{q\,\sin z\over 4c}\,
{2c\over q(1+\cos z)}={\sin w\over 2(1+\cos w)}+{\sin z\over 2(1+\cos z)}\,.
\end{array} \label{6.2}
\eeq
Notice that (\ref{6.2}) defines the values of $u_t(\tau,\cdot)$
at almost every point $x\in\R$, i.e.~at all points outside the support of
the singular part of the measure $\mu_\tau$ defined at (\ref{5.10}).
By the inequality (\ref{1.7}), recalling that $c(u)\geq \kappa^{-1}$,
we obtain
\beq
\int_\R \big|u_t(\tau,x)\big|^2\,dx\leq \kappa^2\,\E(\tau)
\leq \kappa^2\E_0. \label{6.3}
\eeq

To prove (\ref{6.1}), let any $\ve>0$ be given.
There exist finitely many disjoint intervals
$[a_i,\,b_i]\subset \R$, $i=1,\ldots,N$, with the following property.
Call $A_i, B_i$ the points on $\gamma_\tau$ such that
$x(A_i)=a_i$, $x(B_i)=b_i$.  Then one has
\beq
\min\big\{ 1+\cos w(P)\,,~1+\cos z(P)\big\}<2\ve\label{6.4}
\eeq
at every point $P$ on $\gamma_\tau$ contained in one of the arcs
$A_iB_i$, while
\beq
1+\cos w(P)>\ve\,,\qquad\qquad 1+\cos z(P)>\ve\,,\label{6.5}
\eeq
for every point $P$ along $\gamma_\tau$, not contained in any of the arcs
$A_iB_i$.
Call $J\doteq \cup_{1\leq i\leq N} [a_i, b_i]$, $J'=\R\setminus J$,
and notice that, as a function of the original variables, $u=u(t,x)$
is smooth in a neighborhood of the set $\{\tau\}\times J'$.
Using Minkowski's inequality and the differentiability of $u$ on $J'$,
we can write
\beq
\begin{array}{rl}
\lim_{h\to 0}&
{1\over h}\left(\int_{\R} \Big| u(\tau+h,x)-u(\tau,x)-h\,u_t(\tau, x)\Big|^p
dx\right)^{1/p}\\[4mm]
&\leq \lim_{h\to 0}
{1\over h}\left(\int_J \Big| u(\tau+h,x)-u(\tau,x)\Big|^p
dx\right)^{1/p}+\left(\int_J \big|u_t(\tau,x)\big|^p\,dx\right)^{1/p}
\end{array}\label{6.6}
\eeq
We now provide an estimate on the measure of the ``bad" set $J$:
\beq
\begin{array}{rcl}
\meas(J)&=&\int_J dx=\sum_i\int_{A_iB_i} {(1+\cos w)\,p\over 4}\,dX
-{(1+\cos z)\,q\over 4}\,dY\\[4mm]
&\leq &2\ve \sum_i \int_{A_iB_i}
{(1-\cos w)\,p\over 4}\,dX
-{(1-\cos z)\,q\over 4}\,dY\\[4mm]
&\leq& 2\ve  \int_{\gamma_\tau}
{(1-\cos w)\,p\over 4}\,dX
-{(1-\cos z)\,q\over 4}\,dY
\leq 2\ve\,\E_0\,.
\end{array}\label{6.7}
\eeq

Now choose $q=2/(2-p)$ so that
${p\over 2}+{1\over q}=1$. Using H\"older's inequality
with conjugate exponents $2/p$ and $q$, and recalling (\ref{5.14}), we obtain
$$
\begin{array}{rcl}
\int_J \Big| u(\tau+h,x)-u(\tau,x)\Big|^p
dx&\leq &\meas(J)^{1/ q}\cdot
\left(\int_J \Big| u(\tau+h,x)-u(\tau,x)\Big|^2
dx\right)^{p/2}\\[4mm]
&\leq &\big[ 2\ve\,\E_0\big]^{1/q}\cdot
\Big(\big\|u(\tau+h,\cdot)-u(\tau,\cdot)\big\|^2_{\L^2}
\Big)^{p/2}\\[4mm]
&\leq& \big[ 2\ve\,\E_0\big]^{1/q}\cdot
\Big(h^2 \big[ 4(\kappa^3+1)\,\E_0\big]
\Big)^{p/2}.
\end{array}
$$
Therefore,
\beq
\limsup_{h\to 0}
{1\over h}\left(\int_J \Big| u(\tau+h,x)-u(\tau,x)-h\Big|^p
dx\right)^{1/p}\leq [2\ve\,\E_0]^{1/pq}\cdot \big[ 4(\kappa^3+1)\,\E_0
\big]^{1/2}.\label{6.8}
\eeq
In a similar way we estimate
$$
\int_J\big|u_t(\tau,x)\big|^p\,dx\leq
\big[\meas(J)\big]^{1/ q}\cdot
\left(\int_J \Big| u_t(\tau,x)\Big|^2
dx\right)^{p/2},
$$
\beq
\left(\int_J\big|u_t(\tau,x)\big|^p\,dx\right)^{1/p}\leq
\meas(J)^{1/ pq}\cdot \big[\kappa^2\,\E_0\big]^{p/2}.\label{6.9}
\eeq 
Since $\ve>0$ is arbitrary, from (\ref{6.6}), (\ref{6.8}) and (\ref{6.9}) 
we conclude
\beq
\lim_{h\to 0}~
{1\over h}\left(\int_{\R} \Big| u(\tau+h,x)-u(\tau,x)-h\,u_t(\tau, x)\Big|^p
dx\right)^{1/p} =0.
\label{6.10}
\eeq
The proof of continuity of the map $t\mapsto u_t$ is similar.
Fix $\ve>0$. Consider the intervals $[a_i,b_i]$ as before.
Since $u$ is smooth on a neighborhood of $\{\tau\}\times J'$,
it suffices to estimate
$$
\begin{array}{cl}
& \limsup_{h\to 0}\,\int\big|u_t(\tau+h,x)-u_t(\tau,x)\big|^p\,dx \\[4mm]
\leq &\limsup_{h\to 0}
\int_J \big|u_t(\tau+h,x)-u_t(\tau,x)\big|^p\,dx\\[4mm]
\leq & \limsup_{h\to 0} \big[\meas(J)\big]^{1/ q}\cdot
\left(\int_J \Big| u_t(\tau+h,\,x)-u_t(\tau,x)\Big|^2 dx\right)^{p/2} \\[4mm]
\leq & \limsup_{h\to 0}
\big[2\ve\,\E_0\big]^{1/ q}\cdot \Big(\big\|u_t(\tau+h,\cdot)\big\|_{\L^2}
+\big\|u_t(\tau,\cdot)\big\|_{\L^2}\Big)^p\\[4mm]
\leq & \big[2\ve\E_0]^{1/q}\, \big[ 4\E_0\big]^p.
\end{array}
$$
Since $\ve>0$ is arbitrary, this proves continuity.
\v
To extend the result to general initial data, such that 
$(u_0)_x, u_1\in\L^2$, we consider a sequence of smooth initial
data, with $(u^\nu_0)_x, u^\nu_1\in \C^\infty_c$,
with
$u_0^n\to u_0$ uniformly,
$(u_0^n)_x\to (u_0)_x$ almost everywhere and in $\L^2$,
$u_1^n\to u_1$ almost everywhere and in $\L^2$.

The continuity of the function $t\mapsto u_x(t,\cdot)$
as a map with values in $\L^p$, $1\leq p <2$, is proved in
an entirely similar way.

\section{Energy conservation}  

\setcounter{equation}{0}

This section is devoted to the proof of Theorem 3, stating
that, in some sense, the total energy of the solution 
remains constant in time. 

A key tool in our analysis is
the {\it wave interaction potential}, defined as
\beq
\Lambda(t)\doteq (\mu_t^-\otimes\mu_t^+)\big\{ (x,y)\,; ~x>y\big\}\,.
\label{7.1}
\eeq
We recall that $\mu_t^\pm$ are the positive measures defined at
(\ref{5.10'}).
Notice that, if $\mu_t^+, \mu_t^-$ are
absolutely continuous w.r.t.~Lebesgue measure, so that
(\ref{5.9'}) holds, then (\ref{7.1}) 
is equivalent to
$$
\Lambda(t)\doteq \f14 \dint_{x>y} R^2(t,x)\,S^2(t,y)\,dxdy\,.$$
\v
{\bf Lemma 1.} {\it The map $t\mapsto \Lambda(t)$ has locally bounded
variation. Indeed, there exists a one-sided
Lipschitz constant $L_0$ such that
\beq
\Lambda(t)-\Lambda(s) ~\leq L_0\cdot (t-s)\qquad\qquad t>s>0\,.
\eeq
}

To prove the lemma, we first give a formal argument, valid
when the solution $u=u(t,x)$ remains smooth.
We first notice that (\ref{2.4}) 
implies
$${d\over dt} \,(4\Lambda(t))
\leq -\int 2c\, R^2S^2\,dx +\int \big(R^2+S^2\big)\,dx\cdot
\int{c'\over 2c} |R^2S-RS^2|\,dx
$$
$$\leq -2\kappa^{-1}\int  R^2S^2\,dx+ 
4\E_0 \left\|{c'\over 2c}\right\|_{\L^\infty}
\int\big|R^2S-RS^2|\,dx \,,
$$
where $\kappa^{-1}$ is a lower bound for $c(u)$.
For each $\ve>0$ we have
$
|R|\leq \ve^{-1/2}+\ve^{1/2}\,R^2
$.
Choosing $\ve>0$ such that
$$
\kappa^{-1}>4\E_0 \left\|{c'\over 2c}\right\|_{\L^\infty}\cdot 2\sqrt\ve\,,
$$
we thus obtain
$$
{d\over dt}(4 \Lambda(t))\leq -\kappa^{-1}
\int R^2S^2\,dx +{16\,\E_0^2\over\sqrt\ve}
\left\|{c'\over 2c}\right\|_{\L^\infty}\,.
$$
This yields the $\L^1$ estimate
$$
\int_0^\tau \int \big(|R^2S|+|RS^2|\big)\,dxdt =
\O(1)\cdot \big[ \Lambda(0)+ \E_0^2\,\tau\big]=
\O(1)\cdot (1+\tau)\E_0^2,
$$
where $\O(1)$ denotes a quantity whose absolute value
admits a uniform bound, depending only on the function $c=c(u)$ and
not on the particular solution under consideration.
In particular, the map $t\mapsto \Lambda(t)$ has bounded variation
on any bounded interval.
It can be discontinuous, with downward jumps.

To achieve a rigorous proof of Lemma 1, we need to reproduce the 
above argument in terms of the variables $X,Y$.
As a preliminary, we observe that for every $\ve>0$
there exists a constant $\kappa_\ve$
such that
\beq
\begin{array}{rl}
     & \left| \sin z(1-\cos w)-\sin w(1-\cos z)\right| \\[3mm]
\leq & \kappa_\ve\cdot \left(\tan^2\frac{w}{2}+\tan^2\frac{z}{2}
\right)(1+\cos w)(1+\cos z)+\ve (1-\cos w)(1-\cos z)
\end{array}\label{7.3}
\eeq
for every pair of angles $w,z$.

Now fix $0\leq s<t$.  Consider the sets $\Gamma_s,\Gamma_t$ as in (\ref{5.9})
and define
$\Gamma_{st}\doteq \Gamma_t\setminus\Gamma_s$.  
Observing that 
$$dxdt = \frac{pq}{8c}(1+\cos w)(1+\cos z)dXdY\,,$$
we can now write
\beq
\begin{array}{rl}
& \int_s^t\int_{-\infty}^\infty \frac{R^2+S^2}{4}\,dxdt \ = \ (t-s)\E_0 \\[4mm]
= & \dint_{\Gamma_{st}} \frac{1}{4}\left(\tan^2\frac{w}{2}+\tan^2\frac{z}{2}
\right)\cdot \frac{pq}{8c}(1+\cos w)(1+\cos z)\,dXdY\,.
\end{array}\label{7.4}
\eeq
The first identity holds only for smooth solutions, but the second
one is always valid. Recalling (\ref{sour}) and (\ref{5.10'}), and then using 
(\ref{7.3})-(\ref{7.4}), we obtain
$$
\begin{array}{l}
\Lambda(t)-\Lambda(s) \leq -\dint_{\Gamma_{st}} 
\frac{1-\cos w}{8} p\cdot \frac{1-\cos z}{8} q\,dXdY \\[4mm]
\quad +\E_0\cdot \dint_{\Gamma_{st}} \frac{c'}{64 c^2} pq\left[ 
\sin z(1-\cos w)-\sin w(1-\cos z)\right]\,dXdY\\[4mm]
\leq -\frac{1}{64}\dint_{\Gamma_{st}} (1-\cos w)(1-\cos z)\,pq\,dXdY\\[4mm]
\quad+\E_0\cdot \dint_{\Gamma_{st}} \frac{c'}{64 c^2} pq\left[ 
\kappa_\ve\cdot \left(\tan^2\frac{w}{2}+\tan^2\frac{z}{2}
\right)(1+\cos w)(1+\cos z)\right.\\[4mm]
\qquad\left. +\ve (1-\cos w)(1-\cos z)\right]\,dXdY\\[4mm]
\leq \kappa (t-s)\,,
\end{array}
$$
for a suitable constant $\kappa$.  This proves the lemma.

\medskip

To prove Theorem 3, consider the three sets
$$\Omega_1\doteq \Big\{(X,Y)\,;~~w(X,Y)=-\pi\,,\qquad z(X,Y)\not= -\pi
\,,\qquad c'\big(u(X,Y)\big)\not= 0\Big\}\,,$$
$$\Omega_2\doteq \Big\{(X,Y)\,;~~z(X,Y)=-\pi\,,\qquad w(X,Y)\not= -\pi
\,,\qquad c'\big(u(X,Y)\big)\not= 0\Big\}\,,$$
$$\Omega_3\doteq \Big\{(X,Y)\,;~~z(X,Y)=-\pi\,,\qquad w(X,Y)= -\pi
\,,\qquad c'\big(u(X,Y)\big)\not= 0\Big\}\,.$$
From the equations (\ref{2.24}), it follows 
that 
\beq
\meas(\Omega_1)=\meas(\Omega_2)=0\,.
\label{7.5}
\eeq
Indeed, $w_Y\not= 0$ on $\Omega_1$ and $z_X\not= 0$ on $\Omega_2$.

Let $\Omega_3^*$ be the set of Lebesgue points of $\Omega_3$.
We now show that
\beq
\meas\Big(\big\{ t(X,Y)\,;~~(X,Y)\in \Omega_3^*\big\}\Big)=0\,.
\label{7.6}
\eeq
To prove (\ref{7.4}), fix any $P^*=(X^*,Y^*)\in\Omega_3^*$ and let 
$\tau=t(P^*)$. We claim that
\beq
\limsup_{h,k\to 0+} ~{\Lambda(\tau-h)-\Lambda(\tau+k)\over h+k}~=~+\infty
\label{7.7}
\eeq

\begin{figure}[htb]
\centering
\includegraphics[width=0.55\textwidth]{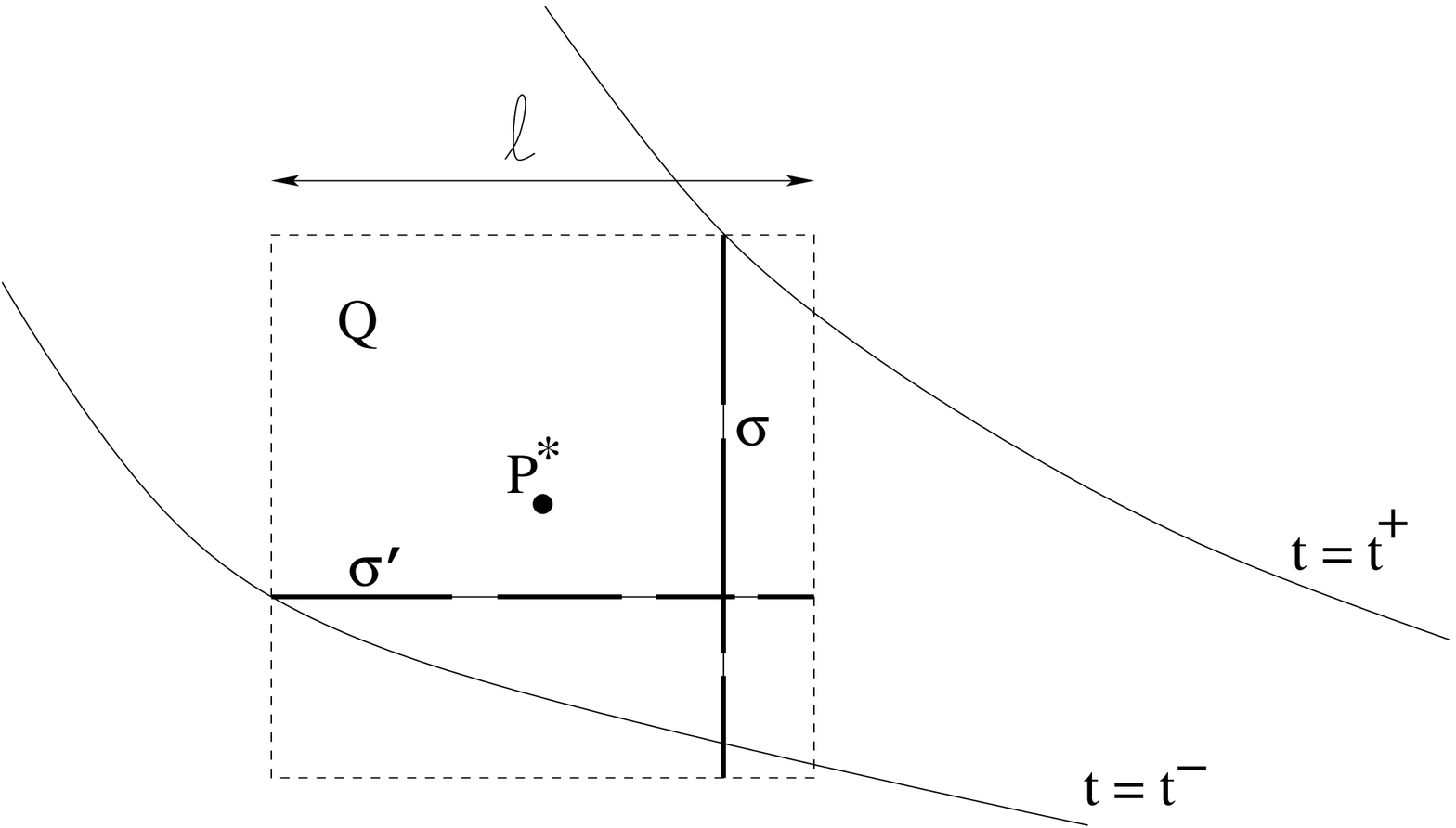} 
\caption{\footnotesize  }
\label{fig6}
\end{figure}

By assumption, for any $\ve>0$ arbitrarily small we can find
$\delta>0$ with the following property.
For any square $Q$ centered at $P^*$ with side of length $\ell<\delta$, 
there exists a
vertical segment $\sigma$ and a horizontal segment $\sigma'$,
as in fig.~6,
such that
\beq 
\meas \big(\Omega_3\cap \sigma\big) \geq (1-\ve)\ell\,,
\qquad\qquad  \meas \big(\Omega_3\cap \sigma'\big) \geq (1-\ve)\ell\,,
\label{7.8}
\eeq
Call 
$$t^+\doteq \max\Big\{ t(X,Y)\,;~~(X,Y)\in \sigma\cup\sigma'\Big\}\,,$$
$$t^-\doteq \min\Big\{ t(X,Y)\,;~~(X,Y)\in \sigma\cup\sigma'\Big\}\,.$$
Notice that, by (\ref{4.2}),
\beq
t^+-t^-\leq \int_\sigma \frac{(1+\cos w)p}{ 4c}\,dX
+\int_{\sigma'} \frac{(1+\cos z)q}{ 4c}\,dY\leq c_0\cdot (\ve\ell)^2\,.
\label{7.9}
\eeq
Indeed, the integrand functions are Lipschitz continuous. Moreover,
they vanish oustide a set of measure $\ve\ell$.
On the other hand, 
\beq
\Lambda(t^-)-\Lambda(t^+)\geq c_1 (1-\ve)^2\ell^2  -
c_2 (t^+-t^-)
\label{7.10}
\eeq
for some constant $c_1>0$. Since $\ve>0$ was arbitrary, 
this implies (\ref{7.5}).

Recalling that the map $t\mapsto \Lambda$ 
has bounded variation, from (\ref{7.5}) it follows  (\ref{7.4}). 

\smallskip
We now observe that the singular part of 
$\mu_\tau$ is nontrivial only if the set 
$$\big\{ P\in\gamma_\tau\,;~~w(P)= -\pi ~~\hbox{or}~~z(P)=-\pi
\big\}$$
has positive 1-dimensional measure.
By the previous analysis, restricted to the region
where $c'\not= 0$, this can happen only for a 
set of times having zero measure.

\vsk

{\bf Acknowledgment:}
Alberto Bressan was supported by the Italian M.I.U.R.,
within the research project \#2002017219, while  
Yuxi Zheng has been partially supported by grants
NSF DMS 0305497 and 0305114.

\vsk

\end{document}